\documentclass[11pt,a4paper]{article}

\usepackage{amsmath, amssymb}
\usepackage{latexsym, color}
\usepackage{pgf,tikz}
\usetikzlibrary{arrows}

\usepackage{epsfig}

\numberwithin{equation}{section}

\begin{document}

\newtheorem{thm}{Theorem}[section]
\newtheorem{cor}[thm]{Corollary}
\newtheorem{lem}[thm]{Lemma}
\newtheorem{prop}[thm]{Proposition}
\newtheorem{defn}[thm]{Definition}
\newtheorem{rem}[thm]{Remark}
\newtheorem{infer}[thm]{Inference}
\newtheorem{Ex}[thm]{EXAMPLE}
\def\nm{\noalign{\medskip}}
\bibliographystyle{plain}


\newcommand{\qed}{\hfill \ensuremath{\square}}
\newcommand{\ds}{\displaystyle}
\newcommand{\pf}{\noindent {\sl Proof}. \ }
\newcommand{\p}{\partial}
\newcommand{\pd}[2]{\frac {\p #1}{\p #2}}
\newcommand{\norm}[1]{\left\| #1 \right\|}
\newcommand{\dbar}{\overline \p}
\newcommand{\eqnref}[1]{(\ref {#1})}
\newcommand{\na}{\nabla}
\newcommand{\one}[1]{#1^{(1)}}
\newcommand{\two}[1]{#1^{(2)}}

\newcommand{\Abb}{\mathbb{A}}
\newcommand{\Cbb}{\mathbb{C}}
\newcommand{\Ibb}{\mathbb{I}}
\newcommand{\Nbb}{\mathbb{N}}
\newcommand{\Kbb}{\mathbb{K}}
\newcommand{\Rbb}{\mathbb{R}}
\newcommand{\Sbb}{\mathbb{S}}

\renewcommand{\div}{\mbox{div}~}

\newcommand{\la}{\langle}
\newcommand{\ra}{\rangle}

\newcommand{\Hcal}{\mathcal{H}}
\newcommand{\Lcal}{\mathcal{L}}
\newcommand{\Kcal}{\mathcal{K}}
\newcommand{\Dcal}{\mathcal{D}}
\newcommand{\Pcal}{\mathcal{P}}
\newcommand{\Qcal}{\mathcal{Q}}
\newcommand{\Scal}{\mathcal{S}}

\def\Ba{{\bf a}}
\def\Bb{{\bf b}}
\def\Bc{{\bf c}}
\def\Bd{{\bf d}}
\def\Be{{\bf e}}
\def\Bf{{\bf f}}
\def\Bg{{\bf g}}
\def\Bh{{\bf h}}
\def\Bi{{\bf i}}
\def\Bj{{\bf j}}
\def\Bk{{\bf k}}
\def\Bl{{\bf l}}
\def\Bm{{\bf m}}
\def\Bn{{\bf n}}
\def\Bo{{\bf o}}
\def\Bp{{\bf p}}
\def\Bq{{\bf q}}
\def\Br{{\bf r}}
\def\Bs{{\bf s}}
\def\Bt{{\bf t}}
\def\Bu{{\bf u}}
\def\Bv{{\bf v}}
\def\Bw{{\bf w}}
\def\Bx{{\bf x}}
\def\By{{\bf y}}
\def\Bz{{\bf z}}
\def\BA{{\bf A}}
\def\BB{{\bf B}}
\def\BC{{\bf C}}
\def\BD{{\bf D}}
\def\BE{{\bf E}}
\def\BF{{\bf F}}
\def\BG{{\bf G}}
\def\BH{{\bf H}}
\def\BI{{\bf I}}
\def\BJ{{\bf J}}
\def\BK{{\bf K}}
\def\BL{{\bf L}}
\def\BM{{\bf M}}
\def\BN{{\bf N}}
\def\BO{{\bf O}}
\def\BP{{\bf P}}
\def\BQ{{\bf Q}}
\def\BR{{\bf R}}
\def\BS{{\bf S}}
\def\BT{{\bf T}}
\def\BU{{\bf U}}
\def\BV{{\bf V}}
\def\BW{{\bf W}}
\def\BX{{\bf X}}
\def\BY{{\bf Y}}
\def\BZ{{\bf Z}}


\newcommand{\Ga}{\alpha}
\newcommand{\Gb}{\beta}
\newcommand{\Gd}{\delta}
\newcommand{\Ge}{\epsilon}
\newcommand{\Gve}{\varepsilon}
\newcommand{\Gf}{\phi}
\newcommand{\Gvf}{\varphi}
\newcommand{\Gg}{\gamma}
\newcommand{\Gc}{\chi}
\newcommand{\Gi}{\iota}
\newcommand{\Gk}{\kappa}
\newcommand{\Gvk}{\varkappa}
\newcommand{\Gl}{\lambda}
\newcommand{\Gn}{\eta}
\newcommand{\Gm}{\mu}
\newcommand{\Gv}{\nu}
\newcommand{\Gp}{\pi}
\newcommand{\Gt}{\theta}
\newcommand{\Gvt}{\vartheta}
\newcommand{\Gr}{\rho}
\newcommand{\Gvr}{\varrho}
\newcommand{\Gs}{\sigma}
\newcommand{\Gvs}{\varsigma}
\newcommand{\Gj}{\tau}
\newcommand{\Gu}{\upsilon}
\newcommand{\Go}{\omega}
\newcommand{\Gx}{\xi}
\newcommand{\Gy}{\psi}
\newcommand{\Gz}{\zeta}
\newcommand{\GD}{\Delta}
\newcommand{\GF}{\Phi}
\newcommand{\GG}{\Gamma}
\newcommand{\GL}{\Lambda}
\newcommand{\GP}{\Pi}
\newcommand{\GT}{\Theta}
\newcommand{\GS}{\Sigma}
\newcommand{\GU}{\Upsilon}
\newcommand{\GO}{\Omega}
\newcommand{\GX}{\Xi}
\newcommand{\GY}{\Psi}

\newcommand{\beq}{\begin{equation}}
\newcommand{\eeq}{\end{equation}}


\title{ An optimal estimate for electric fields on the\\ shortest line segment between \\ two spherical insulators  in three dimensions}

\author{  KiHyun Yun\thanks{Department of Mathematics, Hankuk University of Foreign Studies,
Yongin-si, Gyeonggi-do 449-791, Republic of Korea (kihyun.yun@gmail.com).}}

\maketitle
\abstract{ We consider a gradient estimate for a conductivity problem whose inclusions are two neighboring insulators in three dimensions.  When inclusions with an extreme conductivity (insulators or perfect conductors) are closely located, the gradient can be concentrated in between inclusions and then becomes arbitrarily large as the distance between inclusions approaches zero. The gradient estimate in between insulators in three dimensions has been regarded as a challenging problem, while the optimal blow-up rates in terms of the distance were successfully obtained for the other extreme conductivity problems in two and three dimensions, and are attained on the shortest line segment between inclusions. In this paper, we establish upper and lower bounds of gradients on the shortest line segment between two insulating unit spheres in three dimensions. These bounds present the optimal blow-up rate of gradient on the line segment which is substantially different from the rates in the other problems.}


\section{Introduction}
Let $B_1$ and $B_2$ be bounded and simply connected domains in $\Rbb^d$, $d=2,3$.  We consider the following conductivity problem: for a given
harmonic function $H$ defined in $\Rbb^d$,
 \beq
 \begin{cases}
 \ds \nabla  \cdot \left( (k-1) \chi_{(B_1 \cup B_2)} + 1\right) \nabla u =0 \quad \mbox{in } \Rbb^d, \\
 u(\Bx)- H(\Bx) = O(|\Bx|^{1-d}) \quad \mbox{as } |\Bx| \to
 \infty, 
 \end{cases} \notag
 \eeq where $\chi$ is the characteristic function. Two inclusions $B_1$ and $B_2$ are conductors with conductivity $k  \neq 1$, embedded in the background with conductivity 1. For an extreme conductivity $k = 0$ or $\infty$, a modified model has been used, see \eqref{gov_equation} and \cite {Y}.   Let $\epsilon$ denote the distance  between $B_1$ and $B_2$, $i.e.$,
 \beq \Ge:= \mbox{dist}(B_1, B_2), \notag \eeq  and we assume that  the distance $\Ge$ is small.

\par The problem is to estimate $|\nabla u|$ in the narrow region in between inclusions. This was raised by Babu\u{s}ka in relation to the study of material failure of composites \cite{bab}. In fiber-reinforced composites which consist of stiff fibrous inclusions and the matrix,  high shear stress concentrations can occur in between closely spaced neighboring inclusions. It is important to estimate the shear stress tensor $\nabla u$, while $u$ means the out-of-plane displacement, and the inclusions $B_1$ and $B_2$ are the cross-sections of fibers. Many studies  have been extensively conducted on the gradient estimate due to such practical importance.

\par Successful results have been achieved in all cases except three dimensional insulators which we consider in this paper. Such successful results can be divided into three cases when $k$ stays away from $0$ and $\infty$, when $k $ degenerates to either $0$ (insulating) or $\infty$ (perfectly conducting) in two dimensions, and when $k=\infty$ in three and higher dimensions.  On the other hand, this paper deals with the exceptional case when $k=0$ in three dimensions.  We prove the occurrence of concentration in the narrow region, and also established the optimal blow-up rate for $|\nabla u|$ on the shortest line segment between two insulating unit spheres in terms of $\epsilon$.

\par  We give a brief description of three successful cases mentioned above.  In the first case when $k$ stays away from $0$ and $\infty$, {\it i.e.}, $c_1 < k < c_2$ for some positive constants $c_1$ and $c_2$, it was proved by Li-Vogelius \cite{LV} that $|\nabla u|$ remains bounded regardless of the distance $\Ge$ between inclusions. The boundedness result was extended to elliptic systems by Li-Nirenberg \cite{LN}.

\par In the second case when $k$ is either $0$ or $\infty$ in the two dimensional problem, the gradient $\nabla u$ can become unbounded as the distance $\Ge$ tends to $0$, and the generic blow-up rate of $|\nabla u|$ is $1/\sqrt{\Ge}$  For two circular inclusions, the blow-up rate $1/\sqrt{\Ge}$  was derived by Budiansky-Carrier \cite{BC}, and Kang-Lim $et.al. $ \cite{AKLim,AKLLL} established the precise dependence of $|\nabla u|$ on $\Ge$, radii of disks and $k \in [0,\infty]$.  For inclusions in a sufficiently general class of shapes in $\mathbb{R}^2$, it was shown by  Yun \cite{Y,Y2} that the blow-up rate $1/\sqrt{\Ge}$ is valid at $k = 0$ or $\infty$, see also \cite{LY2} for an enhancement of concentration. Taking it a step further, an asymptotic for the distribution $\nabla u$ was established by Kang-Lim-Yun \cite{KLY}, when $B_1$ and $B_2$ are disks.  For sufficiently general shapes of inclusions, recent results by Ammari $et.al.$ \cite {ACKLY} and Kang-Lee-Yun \cite{KLeeY} yield  a numerically stable method to well describe the asymptotic behavior of $\nabla u$  in $\mathbb{R}^2$.

\par In the third case when $k=\infty$ in three and higher dimensions, Bao-Li-Yin  \cite{BLY} proved that the  generic blow-up rate for the perfectly conducting inclusions is $|\Ge \log \Ge|^{-1}$ in three dimensions and $|\Ge|^{-1}$ in higher ones, see also  \cite{BLY2}  for multiple inclusions. Lim-Yun \cite{LY} also found the explicit dependency of $|\nabla u|$ on the radii as well as the distance $\epsilon$, when two inclusions are spheres in three and higher dimensions, see also \cite{lekner10,lekner11,lekner}.  Morever, an asymptotic for the distribution $\nabla u$ was established by Kang-Lim-Yun \cite{KLY13} and Lim-Yu \cite{LYu}  when $B_1$ and $B_2$ are two perfectly conducting spheres in $\mathbb{R}^3$.     

\par It is important to estimate $|\nabla u|$ on the shortest line segment between inclusions.  In the second and third cases mentioned above, the asymptotic behaviors of $\nabla u$ show in \cite{KLY,ACKLY,KLeeY,KLY13,LYu} that the generic blow-up rates can be attained on the shortest line segment. Since high concentration of $\nabla u$ results from the reflections repeated infinitely between inclusions, refer to \cite{AKLim,Y,LY}, it is reasonable to estimate $|\nabla u|$ on the shortest line segment to obtain the blow-up rate.

\par This paper is mainly concerned with the exceptional case when the inclusions is insulators  ($k=0$) in three dimensions.  The gradient estimate in this case has been regarded as a challenging problem.  An upper bound of $|\nabla u|$ with order $ 1 /{\sqrt {\epsilon}}$ was derived by Bao-Li-Yin \cite{BLY2} in this problem.  To our best knowledge, there has not been any updated or improved result yet.

\par  In this paper,  we prove that the optimal blow-up rate  of  $|\nabla u|$ is $1/ {\epsilon^{\frac {2 -\sqrt 2 }2}} $ on the shortest line segment between two insulating unit spheres in three dimensions. To do so, we establish the upper and lower bounds of $|\nabla u|$ on the  line segment. It is worthy mentioning that the blow-up rate  is substantially different from the rates known in the other extreme conductivity problems.  In terms of methodology, a new method is employed  to derive the optimal blow-up rate due to  the different nature of the problem.

\section{Main Results}
In this paper, we assume that $B_1$ and $B_2$ are a pair of unit spheres $\epsilon$ apart as follows:
 $$B_1 = B_1 \left(-1-\frac {\epsilon} 2 ,0,0\right) \mbox{ and }B_2 = B_1 \left(1+ \frac {\epsilon} 2 ,0,0\right).$$
The quantity $\epsilon$ means the distance between $B_1$ and $B_2$, and is supposed to be small as mentioned in the introduction. The centers of $B_1$ and $B_2$ lie on the $x$ axis.  For any harmonic function $H$ defined on $\mathbb{R}^3$, let $u$ be the solution to the following conductivity problem whose conductivity $k$ degenerates to $0$:
\begin{equation}
\quad \left\{
\begin{array}{ll}
\ds \Delta u  = 0\quad&\mbox{in }\mathbb{R}^3 \setminus \overline {(B_1 \cup B_2)}, \\
\ds \partial_{\nu} u = 0\quad&\mbox{on }  \partial {(B_1 \cup B_2)},\\
\ds u(\Bx) - H (\Bx) = O \left(\frac1 {|\Bx|^2}\right)\quad&\mbox{as } |\Bx| \rightarrow \infty.
\end{array}
\right. \label{gov_equation}
\end{equation} Here, $ \Bx = (x,y,z)$ in $\mathbb{R}^3$.

\par This paper has two main results that provide an optimal blow-up rate $  1 / { \epsilon ^{\frac {2 - \sqrt 2 } {2} }}$ for $|\nabla u |$  on the shortest line segment between two insulating $B_1$ and $B_2$.

\vskip 10pt

\begin{thm} [Upper Bound] \label{main_thm} For any harmonic function $H $ defined in $\mathbb{R}^3$, let $u$ be the solution to \eqref{gov_equation}. There is a constant $C$ independent of the distance $\epsilon >0$ such that 
$$| \nabla  u (x,0,0)|  \leq C \frac 1 { \epsilon ^{\frac {2 - \sqrt 2 } {2} }}$$ for $|x|< \frac {\epsilon} 2$ and any small $\epsilon >0$.
\end{thm}
The proof of this theorem is presented immediately after Proposition \ref{prop_4}, since it is based on Propositions \ref{prop_2}, \ref{prop_3} and \ref{prop_4}.

\vskip 10pt

\begin{thm} [Lower Bound]\label{main_thm_lower} Assume that $$H(x,y,z) = y$$ in $\mathbb{R}^3$. Let $u$ be the solution to \eqref{gov_equation} for $H$. There is a constant $C$ independent of the distance $\epsilon >0$ such that 
$$| \nabla  u (x,0,0)|  \geq C \frac 1 { \epsilon ^{\frac {2 - \sqrt 2 } {2} }}$$ for $|x|< \frac {\epsilon} 2$ and any small $\epsilon >0$.
\end{thm} The proof of Theorem \ref{main_thm_lower} shall be provided in Section \ref{sec5} by modifying the proof of Proposition \ref{prop_4}.  
\vskip 10pt

\par To show Theorem \ref{main_thm},  we use Propositions \ref{prop_2}, \ref{prop_3} and \ref{prop_4}. The first proposition means that the directional derivative  $\partial_x u$ is bounded regardless of $\epsilon$, where  the $x$ axis passes through two centers of $B_1$ and $B_2$.
\begin{prop} \label{prop_2} Let $u$ be the solution to \eqref{gov_equation}  for a harmonic function $H $  in $\mathbb{R}^3$ as given in Theorem \ref{main_thm}.  There is a constant $C$ independent of $\epsilon >0$ such that  $$ \left|\p_{x} u (x,0,0)\right| \leq C $$ for $|x|<\frac \epsilon 2 $ and any small $\epsilon >0$.
\end{prop}
This proof of Proposition \ref{prop_2} is provided in  Subsection  \ref{3-1prop22}.

\vskip 10pt 
\par To estimate the other directional derivatives, we use an upper bound in the second proposition which simplifies this problem.

\begin{prop} \label{prop_3} For a harmonic function $H$  defined in $\mathbb{R}^3$, let $u$ be the solutions to \eqref{gov_equation} for $H$. Then, there exists a large constant $M$ independent of $\epsilon >0$ such that for $ |x|<\frac \epsilon 2 $ and $ \epsilon >0$, $$ \left|\p_{y} u (x,0,0)\right| \leq  \left|\p_{y} u_1 (x,0,0)\right|  ,$$ where $u_1$ is the solutions to \eqref{gov_equation} for $H_1 (x,y,z)= M y$. 
\end{prop}
In  Subsection  \ref{3-2prop23}, we prove Proposition \ref{prop_3}. 
\vskip 10pt 

\par The following proposition is an essential part of this paper that actually yields the first main result.
\begin{prop}\label{prop_4} Let $u$ be the solution to \eqref{gov_equation}  for $H(x,y,z)=y$. Then, there is a constant $C$ independent of $\epsilon >0$ such that 
$$| \nabla  u (x,0,0)|  \leq C \frac 1 { \epsilon ^{\frac {2 - \sqrt 2 } {2} }}$$ for $|x|< \frac {\epsilon} 2$ and any small $\epsilon > 0$.
\end{prop}
In Section \ref{sec4}, we prove Proposition \ref{prop_4}. To do so, we present Proposition \ref{prop_2_2} that obviously implies Proposition \ref{prop_4}. Thus, Section \ref{sec4}  is mainly devoted to showing Proposition \ref{prop_2_2}.

\par Now, we are ready to prove the first main result by virtue of three propositions above.
\vskip10pt

\par 
{\noindent {\sl Proof of Theorem \ref{main_thm} }. \ }  Propositions \ref{prop_2} implies that the directional derivative $\p_ x u$  is bounded independently of $\epsilon>0$ on the line segment between $B_1$ and $B_2$, supposed that the distance $\epsilon$ is small enough. Next, we consider the other directional derivatives. By Propositions \ref{prop_3} and \ref{prop_4}, $$| \p_y  u (x,0,0)|  \leq | \nabla   u_1 (x,0,0)| \leq C_1 \frac 1 { \epsilon ^{\frac {2 - \sqrt 2 } {2} }}$$  for $|x|< \frac {\epsilon} 2$ and any small $ \epsilon>0$,  where 
$u_1$ is the solution to \eqref{gov_equation} for $H=M_1 y$ for a large $M_1>0$. Similarly, we can choose a positive constant $M_2$ independent of $\epsilon>0$ so that $$| \p_z  u (x,0,0)|  \leq | \nabla   u_2 (x,0,0)| \leq C_2 \frac 1 { \epsilon ^{\frac {2 - \sqrt 2 } {2} }}$$  for $|x|< \frac {\epsilon} 2$ and any small $ \epsilon>0$ where 
$u_2$ is the solution to \eqref{gov_equation} for $H=M_2 z$ for a large $M_2>0$.  Therefore, we complete the proof. \qed
\section {Representations of the solution $u$}
Two representations \eqref {sol_decompos_2015} and \eqref {this_lemma_series} of the solution $u$ are introduced in this section. The first representation is used to prove Proposition  \ref{prop_4}, and the second one is for Propositions \ref{prop_2} and \ref{prop_3} whose proofs are presented in Subsections \ref{3-1prop22} and \ref{3-2prop23}, respectively.

\par First, the solution can be  decomposed into three harmonic functions as \beq u (x,y,z)= H (x,y,z) + R_{B_1} (x,y,z) + R_{B_2} (x,y,z), \label{sol_decompos_2015}\eeq where the harmonic function $R_{B_i}$ is defined uniquely in $\mathbb{R}^3 \setminus \overline{B_i}$ and satisfies the decay condition $R_{B_i} = O \left(\frac 1 {|\Bx|^2}\right) $ as $|\Bx| \rightarrow \infty$ for $i=1,2$.  The decomposition can be derived from the representation of $u$ as a sum of $H$ and two single layer potentials on $\p B_1$ and on $\p B_2$, respectively. For details, refer to the invertibility of $ -\frac 1 2 I - \mathbb{K}^* $ presented in Section 2 in \cite{ACKLY}, where $\mathbb{K}^*$ is the  the Neumann-Ponicar{\' e} operator.  It is worthy mentioning that  $R_{B_1}$ is the reflection of  $H + R_{B_2} $ occurring on the insulating inclusion $B_1$ only, and similarly $R_{B_2}$ is the reflection of  $H + R_{B_1} $ on $B_2$.

\par In this paper two harmonic functions $R_{B_1}$ and $R_{B_2}$ play a important rule. They are used for proving Proposition \ref{prop_4} that is actually the first main result in this paper. We study the properties of $R_{B_1}$ and $R_{B_2}$ in Section \ref{sec4} where the proof of Proposition \ref{prop_4} is presented.

\par Second, another representation of $u$ is also introduced in Lemma \ref {another_repre_sol}. This involves the derivations of Propositions \ref{prop_2} and \ref{prop_3}. To illustrate the second representation, we consider  the reflection only for a single inclusion $B_0$ that denotes  the unit sphere with  center $(0,0,0)$, i.e., 
$$ B_0 = B_1 (0,0,0) .$$ For any harmonic function $h$ defined in a neighborhood containing $\overline {B_0}$, let $R_0 (h)$ be the reflection of  $h$ with respect to $B_0$, i.e.,
 \begin{equation}\label{R_0_eqn_g}
\quad \left\{
\begin{array}{ll}
\ds \Delta R_0 (h)  = 0\quad&\mbox{in }\mathbb{R}^3 \setminus \overline {B_0}, \\
\ds \partial_{\nu}  (h+ R_0 (h)) = 0\quad&\mbox{on }  \partial {B_0},\\
\ds R_0 (h) (\Bx)= O \left(\frac1 {|\Bx|^2}\right) &\mbox{as } |\Bx| \rightarrow \infty.
\end{array}
\right.
\end{equation} In the spherical coordinate system, 
\begin {equation} R_0( h ) (\rho, \theta, \phi) = \frac 1 {\rho} h \left(\frac 1 {\rho}, \theta, \phi\right) - \int_{0} ^{\frac 1 {\rho}  }  h \left(s, \theta, \phi\right) ds \label {eq_R_0_1}\end{equation}
 for $\rho \geq 1$, where  $(x,y,z)=\left( \rho \cos  \theta \sin \phi, \rho \sin  \theta \sin \phi, \rho  \cos \phi  \right)$ in $\mathbb{R}^3$. In the Cartesian coordinate system, 
 \beq \p_{y} R_0( h ) (x, 0, 0) = \frac 1 {x^3} \p_{y} h \left(\frac 1 x,0,0\right) - \frac 1 x \int_{0} ^{\frac 1 x  }  s \p_{y} h \left(s,0,0\right) ds \label{R_0_deriv_1} \eeq
for $x \geq 1$.

\par Similarly, for any harmonic function $h$ defined in a neighborhood containing $\overline {B_i}$,  we define $R_i (h)$ as the reflection of $h$ with respect to $B_i$ for $i=1,2$ as follows:
 \begin{equation}
\quad \left\{
\begin{array}{ll}
\ds \Delta R_i (h)  = 0\quad&\mbox{in }\mathbb{R}^3 \setminus \overline {B_i},\notag \\
\ds \partial_{\nu}  (h+ R_i (h)) = 0\quad&\mbox{on }  \partial {B_i},\notag\\
\ds R_i (h) (\Bx)= O \left(\frac1 {|\Bx|^2}\right) &\mbox{as } |\Bx| \rightarrow \infty.\notag
\end{array}
\right. \label{3-5_2015}
\end{equation}

\begin{lem}\label{another_repre_sol}  The following sum converges to the solution $u$ in the sense of the Sovolev space $W^{4,\infty} (\mathbb{R}^3 \setminus \overline{(B_1 \cup B_2)}) $, i.e., \begin{align} u(\Bx)=& H (\Bx)+ R_1(H) (\Bx) + R_2(H) (\Bx) \notag\\ &+ \sum_{n=1} ^{\infty} (R_1 R_2)^n (H) (\Bx)  + (R_1 R_2)^n R_1(H) (\Bx) \notag \\ &+ \sum_{n=1} ^{\infty}  (R_2 R_1)^n (H) (\Bx)  + (R_2 R_1)^n R_2(H) (\Bx)  \label{this_lemma_series}\end{align} for any $\Bx \in \mathbb{R}^3 \setminus \overline{(B_1 \cup B_2)} $.  \end{lem}
\pf   We begin by showing two properties \eqref{extra3-6} and \eqref{2nd_stan_r} that are essential to prove the convergence of the series \eqref{this_lemma_series}. Let $h$ be a harmonic function defined in $H^1\left(\mathbb{R}^3 \setminus {B_1}\right)$ with the decay condition  $h(\Bx) = O(\frac 1{|\Bx|^2})$ as $|\Bx| \rightarrow  \infty$.

\par We first show that \beq  \left| \int_{\p B_2} R_2 (h ) \p_{\nu} \left( R_2 (h) \right)~dS\right| \leq \frac 1 {(1+ \epsilon)^3  }  \left| \int_{\p B_1} h \p_{\nu} h ~dS\right|.\label {extra3-6}\eeq  Since ${B_1}$ is the unit sphere with the center  $  \left(- 1 -  \frac  \epsilon 2, 0,0 \right)$,  the function $h$ can be expressed in terms of spherical harmonic functions whose  center  is $  \left(- 1 -  \frac \epsilon 2, 0,0 \right)$. By the decay condition of $h$,  we have 
\begin{align*}
\frac 1 {(1+\epsilon)^3}    \left| \int_{\p B_1} h \p_{\nu} h ~dS\right|& \geq  \left| \int_{|\Bx - (- 1 - \frac 1 2 \epsilon,0,0 )  | = 1+ \epsilon  } h \p_{\nu} h ~dS \right| \\& =  \left|\int_{|\Bx - (- 1 - \frac 1 2 \epsilon, 0,0 )  | > 1+ \epsilon  } |\nabla h |^2 ~dV \right|\\& \geq  \left| \int_{\p B_2   } h \p_{\nu}h ~dS \right| \\& \geq  \left| \int_{\p B_2   } R_2 (h) \left( \p_{\nu} R_2 (h) \right)  ~dS \right|.
\end{align*} Thus, we have \eqref {extra3-6}. 

\par Second, we prove the existence of a constant $C_1$ such that  \beq
 C_1 \left| \int_{\p B_1} h \p_{\nu} h ~dS\right|^{\frac 1 2 }  \geq \norm{R_2 (h)}_{{W}^{4,\infty}\left(\mathbb{R}\setminus B_2\right)}. \label {2nd_stan_r}\eeq The mean value property for harmonic functions yields the inequality  \begin{align*}
 \left| \int_{\p B_1} h \p_{\nu} h ~dS\right|^{\frac 1 2 }&  =  \left|\int_{ \mathbb {R}^3 \setminus B_1 } |\nabla h |^2 ~dV \right|^{\frac 1 2 } \\& \geq  C_2 \norm { \nabla h }_{L^{\infty} ({|\Bx - ( 1 + \frac 1 2 \epsilon,0,0 )  | \leq  1+ \frac 1 2 \epsilon  })} \\ & \geq 
C_3 \left( \max_{|\Bx - ( 1 + \frac 1 2 \epsilon,0,0 )  | \leq  1+ \frac 1 2 \epsilon  } ( h (\Bx)) - \min_{|\Bx - ( 1 + \frac 1 2 \epsilon,0,0 )  | \leq  1+ \frac 1 2 \epsilon  }  ( h (\Bx)) \right).
\end{align*} The positive constants $C_2$ and $C_3$ are used above regardless of  choosing a harmonic function $h$  in $H^1 (\mathbb{R}^3 \setminus \overline{B_1})$. Note that $R_2(h)$ can be extended as a harmonic function into the area $|\Bx - ( 1 + \frac 1 2 \epsilon,0,0 )  | >  \frac 1{1+\frac 1 2 \epsilon} $, and that the analogous formula for $R_2(h)$ with \eqref {eq_R_0_1} is valid in the extended domain. Thus, the  formula for $R_2(h)$ implies 
$$
 C_4 \left| \int_{\p B_1} h \p_{\nu} h ~dS\right|^{\frac 1 2 } \geq \max_{|\Bx - ( 1 + \frac 1 2 \epsilon,0,0 )  | \geq  \frac 1 {1+ \frac 1 2 \epsilon}  }|R_2 (h) (\Bx)|.$$
The distance between the boundaries of  $B_2 $ and the extended domain of $R_2(h)$ is greater than $\frac 1 4 \epsilon$. For any point in $\mathbb{R}^3 \setminus B_2$,  the harmonic function  $R_2(h)$ is defined in the open sphere whose center is the given point and the radius is $\frac 1 4 \epsilon$.
Thus, a gradient estimate for harmonic functions implies $$ C_5 \left| \int_{\p B_1} h \p_{\nu} h ~dS\right|^{\frac 1 2 }  \geq \norm{\nabla R_2 (h)}_{L^{\infty}\left(\mathbb{R}^3 \setminus B_2\right)},$$ and moreover, the bound \eqref {2nd_stan_r} for the higher order derivatives can also be derived in the same way.

\par Now, we are ready to prove this lemma.  Applying \eqref {extra3-6} repeatedly, we have 
\begin{align*} & \left|\int_{\p B_1}   (R_1 R_2)^n (H)   \p_{\nu} \left ( (R_1 R_2)^n (H) \right) dS \right| \\& +  \left|\int_{\p B_1}   (R_1 R_2)^n (R_1 H)   \p_{\nu} \left ( (R_1 R_2)^n (R_1 H) \right) dS \right|  \\  & + \left|\int_{\p B_2}   (R_2 R_1)^n (H)   \p_{\nu} \left ( (R_2 R_1)^n (H) \right) dS \right| \\& +  \left|\int_{\p B_2}   (R_2 R_1)^n (R_2 H)   \p_{\nu} \left ( (R_2 R_1)^n (R_2 H) \right) dS \right| \\& \leq \frac {2}{(1+\epsilon)^{3(2n-1)}}  \sum_{i=1}^2  \left|\int_{\p B_i}   R_i  (H)   \p_{\nu} R_i (H)  dS \right| \\& \leq \frac {2}{(1+\epsilon)^{3(2n-1)}} \sum_{i=1}^2  \left|\int_{\p B_i}  H   \p_{\nu} H  dS \right|  \end{align*} for any $n=1,2,3,\cdots$. By \eqref{2nd_stan_r}, we have \begin{align} &\norm{ (R_1 R_2)^n (H) }_{{W}^{4,\infty}\left(\mathbb{R}^3\setminus B_1\right)} +  \norm{ (R_1 R_2)^n (R_1 H)  }_{{W}^{4,\infty}\left(\mathbb{R}^3\setminus B_1\right)} \notag  \\&+ \norm{ (R_2 R_1)^n (H) }_{{W}^{4,\infty}\left(\mathbb{R}^3\setminus B_2\right)} +  \norm{ (R_2 R_1)^n (R_2 H)  }_{{W}^{4,\infty}\left(\mathbb{R}^3\setminus B_2\right)}   \notag  \\ &\leq C_1 \frac {2}{(1+\epsilon)^{3(2n-1)}} \sum_{i=1}^2  \left|\int_{\p B_i}  H   \p_{\nu} H  dS \right|  \label{abcdefghi} \end{align}  for any $n=1,2,3,\cdots$, where $C_1$ is the constant in \eqref {2nd_stan_r}. This implies that the series in the right hand side of \eqref{this_lemma_series} are convergent in the sense of $W^{4,\infty} (\mathbb{R}^3 \setminus (B_1 \cup B_2)) $ and thus satisfies \eqref{gov_equation}. Hence, the series converges to the solution $u$.\qed

\vskip 10pt
\par Propositions  \ref{prop_2} and  \ref{prop_3} can be derived from basic properties of the representations introduced in this section.

\subsection {Proof of Proposition \ref{prop_2}} \label{3-1prop22}

We begin in considering the case of a single inclusion $B_0$.  As defined before, $R_0(h)$ is the reflection of a given harmonic function $h$ with respect to $B_0$.  By the equation \eqref {eq_R_0_1}, 
\beq \p_x R_0 (h) (x,0,0) = -\frac 1 {x^3} \p_x h \left(\frac 1 x , 0 ,0\right)\mbox{ for }x \geq 1.\label{pxr0}\eeq Here, $\frac 1 x$ in the equation means the first coordinate of  $\left(\frac 1 x,0,0\right)$ at which the image charge of $(x,0,0)$ with respect to $B_0$ is located.

\par  Lemma \ref{another_repre_sol} implies that  the solution $u$ results from the recursive reflections on $B_1$ and $B_2$. Dealing with the recursive reflections, we define $r_1(x)$ and $r_2 (x)$ as the first coordinates of the image charges of $(x,0,0)$ with respect to $B_1$ and $B_2$, respectively. Thus, 
\begin {align*} 
r_1 (x)&= \frac 1 {x + 1 + \frac \epsilon 2 } - \left(1 + \frac \epsilon 2 \right) ~&&\mbox{ for }x \geq -\frac \epsilon 2 , \\r_2 (x)& = 1 + \frac \epsilon 2  -\frac 1 { 1 + \frac \epsilon 2 - x } ~&&\mbox{ for }x \leq \frac \epsilon 2 .
\end{align*}  For $|x|<\frac \epsilon 2 $, we define two sequaneces $(r_{A n})$ and $(r_{B n})$ as
$$ r_{A 2n-1} = (r_2 r_1)^{n-1}r_2 (x)~\mbox{ and }~ r_{A2n} = -(r_1 r_2)^{n}(x), $$
and 
$$ r_{B 2n-1} = -(r_1 r_2)^{n-1}r_1 (x)~\mbox{ and }~ r_{B2n} = (r_2 r_1)^{n}(x) $$
for $n=1,2,3,\cdots$, where $(r_i r_j) (x) = r_i (r_j(x))$ for $\{i,j\}=\{1,2\}$.  Applying \eqref {pxr0} to \eqref{this_lemma_series},  \begin{align*} \p_x u(x,0,0)=& \p_x H (x,0,0)\\&+ \sum_{n=1}^{\infty} (-1)^n \left( \prod_{k=1}^n \left(1+\frac \epsilon 2 - r_{Ak}\right) \right)^3 \p_x H\left((-1)^{n+1}r_{An},0,0\right)\\&+ \sum_{n=1}^{\infty} (-1)^n \left( \prod_{k=1}^n \left(1+\frac \epsilon 2 - r_{Bk}\right) \right)^3 \p_x H\left((-1)^n r_{Bn},0,0\right) \end{align*} for any $|x|\leq \frac \epsilon 2$.

 \par  Indeed, two positive sequences $(r_{An})$ 
 and $(r_{Bn})$ are increasing and converge to a number that is $\sqrt \epsilon + O(\epsilon)$. There are some properties that can be shown easily:  
$$ \frac  \epsilon 2  \leq r_{An} (x) \leq 2\sqrt {\epsilon},~~~\frac  \epsilon 2 \leq r_{Bn} (x) \leq 2\sqrt {\epsilon} $$
 for $|x|\leq \frac \epsilon 2 $ and $n=0,1,2,3,\cdots$, and  
$$ \frac 1 {40 } \sqrt \epsilon \leq  (r_2 r_1)^n (x)  ~\mbox{~and~}~ \frac 1 {40 } \sqrt \epsilon \leq - (r_1 r_2)^n  (x)  $$ for any $n >\frac 1 {20\sqrt {\epsilon}}$ and $|x|\leq \frac {\epsilon}{2}$.  These properties follow immediately from Lemma \ref{prop_x_n} in this paper.

\par Now, we prove the boundedness in Proposition \ref{prop_2}. Let $\alpha_A = -1 $ and $\alpha_B=1$.  For $j=A$ or $B$,
\begin{align*}
&\left|\sum_{n=1}^{\infty} (-1)^n \left( \prod_{k=1}^n \left(1+\frac \epsilon 2 - r_{jk}\right) \right)^3 \p_x H(\alpha_j (-1)^{n} r_{jn},0,0)\right|\\&\leq \sum_{{\tilde n}=1}^{\infty} \left( \prod_{k=1}^{2{\tilde n}-1} \left(1+\frac \epsilon 2 - r_{jk}\right) \right)^3 \\&~~~~~~~~ \times\Big(\left| \p_x H(-\alpha_j r_{j(2{\tilde n}-1)},0,0) - \p_x H(\alpha_j r_{j (2\tilde n)},0,0) \right| + (6\sqrt \epsilon + O(\epsilon))\big| \p_x H(r_{j (2\tilde n)},0,0)\big|  \Big)\\&\leq \sum_{{\tilde n}=1 } ^{ \left[\frac 1 {20 \sqrt {\epsilon}}\right]} \left (4\sqrt \epsilon \norm {\p_x ^2 H (x,0,0)}_{L^{\infty} ([-2 \sqrt \epsilon ,2 \sqrt \epsilon  ])} +  (6\sqrt \epsilon + O(\epsilon)) \norm {\p_x H (x,0,0)}_{L^{\infty} ([-2 \sqrt \epsilon ,2 \sqrt \epsilon  ])}\right)\\&~~+\sum_{{\tilde n}= \left[\frac 1 {20 \sqrt {\epsilon}}\right]+1}^{\infty} 
8 \left(1 - \frac {\sqrt \epsilon}{40} +\frac \epsilon 2\right)^{6\left({\tilde n}- \left[\frac 1 {20 \sqrt {\epsilon}}\right] \right)}\\&~~~~~~~~~~~ \times\left (4\sqrt \epsilon \norm {\p_x ^2 H (x,0,0)}_{L^{\infty} ([-2 \sqrt \epsilon ,2 \sqrt \epsilon  ])} +  (6\sqrt \epsilon + O(\epsilon)) \norm {\p_x H (x,0,0)}_{L^{\infty} ([-2 \sqrt \epsilon ,2 \sqrt \epsilon  ])}\right) \\&\leq C \left ( \norm {\p_x ^2 H (x,0,0)}_{L^{\infty} ([-2 \sqrt \epsilon ,2 \sqrt \epsilon  ])} +  \norm {\p_x H (x,0,0)}_{L^{\infty} ([-2 \sqrt \epsilon ,2 \sqrt \epsilon  ])}\right).  \end{align*} Hence, we have the boundedness of $|\p_x u (x,0,0)|$ as 
$$\norm {\p_x u (x,0,0)}_{L^{\infty} ([-\frac \epsilon 2, \frac  \epsilon  2])} \leq  C \left ( \norm {\p_x ^2 H (x,0,0)}_{L^{\infty} ([-2 \sqrt \epsilon ,2 \sqrt \epsilon  ])} +  \norm {\p_x H (x,0,0)}_{L^{\infty} ([-2 \sqrt \epsilon ,2 \sqrt \epsilon  ])}\right) $$ for small $\epsilon>0$. \qed

\subsection {Proof of Proposition \ref{prop_3}} \label{3-2prop23}
In the following lemma, we first consider the model of a single inclusion $B_0$ that is much simpler than our model of two inclusions. Second, applying the lemma to Lemma \ref{another_repre_sol}, we prove Proposition \ref{prop_3}.
\begin{lem}\label{lem_prop_3} Let $h$ be a harmonic function defined in a neighborhood containing the closure of the unit sphere $ {B_0}= B_1 (0,0,0)$, and let $R_0 (h)$ be defined as \eqref{R_0_eqn_g}. Suppose that 
$$ \p_y h (x,0,0) \geq 0 ~\mbox{ and }~ \p_x \p_y h (x,0,0)\geq 0 $$ for $ 0\leq x \leq 1$. Then, 
$$ \p_y R_0( h) (x,0,0) \geq 0  ~\mbox{ and }~ \p_x \p_y R_0( h )(x,0,0) \leq 0 $$  for $ x \geq 1$.
\end{lem} \pf   By \eqref{R_0_deriv_1} and the assumption,  \begin{align*}\p_{y} R_0( h ) (x, 0, 0) &= \frac 1 {x^3} \p_{y} h \left(\frac 1 x,0,0\right) - \frac 1 x \int_{0} ^{\frac 1 x  }  s \p_{y} h \left(s,0,0\right) ds \notag \\&= \frac 1 x \int_{0} ^{\frac 1 x  }  s  \left( 2 \p_{y} h \left(\frac 1 x,0,0\right) - \p_{y} h \left(s,0,0\right)\right) ds  \geq 0 \end{align*} for $ x \geq 1$. Thus, we have the first bound.

\par We can get the second bound from \eqref{R_0_deriv_1} as follows:  \begin{align*} 
\p_{x}\p_{y} R_0( h) (x, 0, 0) = &  - \frac 2 {x^4}   \p_{y} h \left(\frac 1 x,0,0\right) + \frac 1 {x^2} \int_{0} ^{\frac 1 x  }  s \p_{y} h \left(s,0,0\right) ds - \frac 1 {x^5}  \p_x  \p_{y} h \left(\frac 1 x,0,0\right)  \notag \\&\leq - \frac 1 {x^2} \int_{0} ^{\frac 1 x  }  s  \left( 4 \p_{y} h \left(\frac 1 x,0,0\right) - \p_{y} h \left(s,0,0\right)\right) ds  \leq  0. \end{align*} \qed
\vskip 10pt 

Now, we are ready to prove Proposition \ref{prop_3}.

\vskip 10pt 
\par {\noindent {\sl Proof of Proposition \ref{prop_3}}. \ } From the definitions of $R_1$ and $R_2$, 
\begin{align}
&\p_y R_i (y) (x,0,0) = \frac 1 {2 \left( 1 + \frac \epsilon 2 - (-1)^i x\right)^3}>0,\label{prop32-proof-2015}\\
&(-1)^i \p_x \p_y R_i (y) (x,0,0) = \frac {3} {2 \left( 1 + \frac \epsilon 2 - (-1)^i x\right)^4}>0\label{prop32-proof-2015_2}
\end{align}for $-1-\frac \epsilon 2 \leq (-1)^i x \leq \frac \epsilon 2$ and $i = 1, ~2$. There is a large $M >0$ such that 
$$ \p_y R_i (H) (x,0,0) \leq \p_y R_i (My) (x,0,0) $$
and 
$$  (-1)^i \p_x \p_y R_i (H) (x,0,0) \leq  (-1)^i \p_x \p_y R_i (My) (x,0,0) $$ for $-1-\frac \epsilon 2 \leq (-1)^i x \leq \frac \epsilon 2$ and $i = 1, 2$, since $R_i (My) = M R_i (y)$.  By mathematical induction, Lemma \ref{lem_prop_3} allows  the upper and lower bounds of  $\p_y (R_1 R_2)^n (H) (x,0,0)$ and $\p_x \p_y (R_1 R_2)^n (H) (x,0,0)$ and so on such that for any $n=1,2,3,4,\cdots$, 
$$  \p_y (R_1 R_2)^n (H) (x,0,0) \leq \p_y (R_1 R_2) (My) (x,0,0) ,$$
$$ \p_x \p_y (R_1 R_2)^n (H) (x,0,0) \geq \p_x \p_y (R_1 R_2)^{n} (My) (x,0,0)  $$ and  $$  \p_y (R_1 R_2)^{n-1} R_1 (H) (x,0,0) \leq \p_y (R_1 R_2)^{n-1} R_1 (My) (x,0,0), $$ $$  \p_x \p_y (R_1 R_2)^{n-1} R_1 (H) (x,0,0) \geq  \p_x \p_y (R_1 R_2)^{n-1} R_1 (My) (x,0,0) $$ for $-\frac \epsilon 2 \leq  x \leq 1+\frac \epsilon 2 $, and $$ \p_y (R_2 R_1)^n (H) (x,0,0) \leq \p_y (R_2 R_1)^n (My) (x,0,0) ,$$ $$\p_x \p_y (R_2 R_1)^n (H) (x,0,0) \leq \p_x\p_y (R_2 R_1)^n (My) (x,0,0) $$ and  $$ \p_y (R_2 R_1)^{n-1} R_1 (H) (x,0,0) \leq \p_y (R_2 R_1)^{n-1} R_2 (My) (x,0,0) ,$$ $$\p_x \p_y (R_2 R_1)^{n-1} R_1 (H) (x,0,0) \leq \p_x \p_y (R_2 R_1)^{n-1} R_2 (My) (x,0,0) $$   for $-1-\frac \epsilon 2 \leq  x \leq \frac \epsilon 2 $.
By Lemma \ref{another_repre_sol}, we have  the upper bound of $\p_y u(x,0,0)$ as
\begin{align*}  \p_y u(x,0,0) \leq & M+ \p_y R_1(My)(x,0,0) +\p_y R_2(My)(x,0,0) \notag\\ &+ \sum_{n=1} ^{\infty} \p_y (R_1 R_2)^n (My)(x,0,0)  + \p_y (R_1 R_2)^n R_1(My)(x,0,0) \notag \\ &+ \sum_{n=1} ^{\infty}  \p_y (R_2 R_1)^n (My)(x,0,0)  + \p_y (R_2 R_1)^n R_2(My)(x,0,0) \notag\\ =&\p_y u_1(x,0,0)\end{align*} for $|x|\leq \frac \epsilon 2$, when $u_1$ is the solutions to \eqref{gov_equation} for $H_1 (x,y,z)= M y$.  The lower bound is also obtained in the same way.

\qed

\section{ Proof of Proposition \ref{prop_4}} \label{sec4}
In this section, we assume that $$H(x,y,z)= y \mbox{ in }\mathbb{R}^3$$ and $u$ is the solution to \eqref{gov_equation}  for $H=y$. As defined in the decomposition \eqref{sol_decompos_2015}, two harmonic functions $R_{B_1}$ and $R_{B_2}$ satisfy \beq u (x,y,z)= y + R_{B_1} (x,y,z) + R_{B_2} (x,y,z), \notag\eeq where $R_{B_i}$ is defined in $\mathbb{R}^2 \setminus \overline{B_i}$ and satisfies the decay condition $R_{B_i} = O \left(\frac 1 {|\Bx|^2}\right) $ as $|\Bx| \rightarrow \infty$ for $i=1,2$.

 \par This section is mainly devoted to proving Proposition \ref{prop_2_2} that obviously implies Proposition \ref{prop_4}, since  $\p_z u (x,0,0) = 0 $  and $\p_x u (x,0,0)$ is bounded for $|x|\leq \frac {\epsilon} 2$ by Proposition \ref{prop_2}.  
 
\begin{prop} \label{prop_2_2}$\p_{y} R_{B_1} \left(x ,0,0\right)$ is decreasing in $[-\frac {\epsilon} 2 , \infty)$ and
\beq 0\leq  \p_{y} R_{B_1} \left(-\frac {\epsilon} 2,0,0\right)\lesssim \frac 1 { \epsilon ^{\frac {2 - \sqrt 2 } {2} }}. \label{prop41thefirst2015} \eeq   It follows from the property $\p_{y} R_{B_1} \left(x ,0,0\right) = \p_{y} R_{B_2} \left(-x ,0,0\right) $ that  
\beq \p_{y} u  \left(x ,0,0\right)  =  1+ \p_{y} R_{B_1} \left(x ,0,0\right) + \p_{y} R_{B_2} \left(x ,0,0\right) \lesssim \frac 1 { \epsilon ^{\frac {2 - \sqrt 2 } {2} }} \label{prop41thefirst2015-2} \eeq  for $|x|<\frac \epsilon 2 $.
\end{prop}    Here and throughout this paper, $a_1 \lesssim b_1 $ means $ a_1 \leq C_1 b_1$ and $a_2 \simeq b_2$ stands for $\frac 1 {C_2} a_2 \leq b_2 \leq C_2 a_2$ for some constants $C_1$ and $C_2$ independent of $\Ge$.

\par The proof of Proposition \ref{prop_2_2} is presented in Subsection \ref{subsection3-2},  based on the lemmas in Subsection \ref{subsection3-1}.

\subsection{ Basic Properties of  $\p_{y} R_{B_1} (x,0,0)$} \label{subsection3-1}

We consider the behavior of $\p_{y} R_{B_1} (x,0,0)$ to derive Proposition \ref{prop_2_2}. In this subsection,  $H(x,y,z)= y \mbox{ in }\mathbb{R}^3$ as assumed early in Section \ref{sec4}. For convenience, we define the function $P:[1,\infty) \rightarrow \mathbb{R}$  as
$$P(x) = \p_{y} R_{B_1} \left(x-1 -\frac {\epsilon} 2 ,0,0\right) $$ that is a horizontal shift of $\p_{y} R_{B_1}(x,0,0)$.  The translation moves the left inclusion $B_1$ to $B_1 (0,0,0)$ so that the domain of  $P$ is  the interval $[1,\infty)$ with the initial point $1$. The symmetry between $\p_{y} R_{B_1} (x,0,0) $ and $ \p_{y} R_{B_2} (x,0,0)$ yields \begin{align*}
\p_{y} u (x,0,0)&= \p_{y} H (x,0,0) + \p_{y} R_{B_1} (x,0,0) + \p_{y} R_{B_2} (x,0,0)\\&= 1+  P \left(x+1+\frac \epsilon 2\right) + P \left(-x+1+\frac \epsilon 2\right)
\end{align*} for $|x|\leq \frac \epsilon 2$, since $H(x,y,z)= H(-x,y,z)$.  We study the behavior of $P(x)$ especially for small $x-1 \geq 0$ to prove the bound \beq P (x)\lesssim \frac 1 { \epsilon ^{\frac {2 - \sqrt 2 } {2} }} \label{main_esti_P}\eeq that means Propositions \ref{prop_2_2} and \ref{prop_4}.

\par The basic properties of $P$ are introduced in this subsection. First, Lemma \ref{Ga} provides a fundamental equation \eqref{B} of $P$ that yields almost all properties of $P$ including the main result in this paper.  Second,  Lemmas \ref{Na} and \ref{Da} describe the geometric behavior of $P$. Third, Lemma \ref{integral} presents an estimate for the integral value of $P$ that determines the blow-up rate of $P$ as $\epsilon$ approachs $0$. Finally,  based on these properties, our main result \eqref {main_esti_P} can be obtained in Subsection \ref{subsection3-2} to show Proposition \ref{prop_2_2}. 

 \par Dealing with the lemmas, it is necessary to define two special points $p_1$ and $p_2$ as a pair of solutions to \beq 2+\epsilon -  \frac 1 x = x ,\label{p_1p_2}\eeq where  $p_1 < p_2$.  Indeed, they are the fixed points of the composition of two Kelvin transforms with respect to $B_1(0,0,0)$ and $B_1 (2+\epsilon, 0,0)$, and can be calculated directly as \beq p_i = 1 + (-1)^i \sqrt \epsilon + O_i( \epsilon) \label{p_1p_2}\eeq
for $i=1,2$.

\par In the following lemma,  we establish a fundamental equation of $P$ which is essential in deriving many properties of $P$ in this paper.
\begin{lem} \label{Ga} The function $P:[1,\infty) \rightarrow \mathbb{R}$  satisfies
\beq \frac 1 {x^3} P\left(2+\epsilon  - \frac 1 x \right)- \frac 1 x \int_{0} ^{\frac 1 x} s P(2+\epsilon - s) ds + \frac 1 2 \frac 1 {x^3} = P(x) \label{B}\eeq  for any $x \geq 1$.
\end{lem}
 \pf  We first  recall the property   \eqref{R_0_deriv_1}  of  the reflection with respect to a single inclusion. Second, we apply  \eqref{R_0_deriv_1}  to the case of two neighboring inclusions thus to derive \eqref{B}.  
 
 \par We first consider  the reflection only for a single inclusion $B_0$ that is the unit sphere $B_1 (0,0,0)$. For any harmonic function $h$ defined in a neighborhood containing $\overline {B_0}$,  the reflection $R_0 (h)$ with respect to $B_0$ satisfies 
 \beq \p_{y} R_0( h ) (x, 0, 0) = \frac 1 {x^3} \p_{y} h \left(\frac 1 x,0,0\right) - \frac 1 x \int_{0} ^{\frac 1 x  }  s \p_{y} h \left(s,0,0\right) ds \label{R_0_deriv} \eeq
for $x \geq 1$, that is actually \eqref{R_0_deriv_1}.

\par Second, we consider the solution $u$ to \eqref{gov_equation}. It can be  decomposed into three harmonic functions as $ u (x,y,z)= H (x,y,z) + R_{B_1} (x,y,z) + R_{B_2} (x,y,z)$. From definition, $ R_{B_1} (x,y,z)$ can be regarded as the reflection of $ h= H (x,y,z) +  R_{B_2} (x,y,z) $ with respect to $B_1 = B_0 +\left(-1 -\frac \epsilon 2 , 0, 0 \right)$, and   
$ R_{B_2} (x,y,z) = R_{B_1} (-x,y,z) $ due to the symmetric property of $H (x,y,z)= y$. Since $P(x) = \p_{y} R_{B_1} \left(x-1 -\frac {\epsilon} 2,0,0\right)$,  the equality \eqref{R_0_deriv} thus yields
\begin{align}
&P(x) = \p_{y} R_{B_1} \left(x-1 -\frac {\epsilon} 2,0,0\right) \notag \\&= \p_{y} R_0\left(  H \left(x -1-\frac {\epsilon} 2,y,z\right )  + R_{B_2} \left(x -1-\frac {\epsilon} 2,y,z\right )\right ) (x, 0, 0)   \notag\\&= \p_{y} R_0\left(  y  + R_{B_1} \left(-x +1 + \frac {\epsilon} 2,y,z\right )\right ) (x, 0, 0)   \notag\\&=\frac 1 {x^3} + \frac 1 {x^3} \p_{y} R_{B_1} \left(-\frac 1 x +1 + \frac {\epsilon} 2,0,0\right ) - \frac 1 {2 x^3}-\frac 1 x \int_{0} ^{\frac 1 x  }  s \p_{y} R_{B_1} \left(-s +1 + \frac {\epsilon} 2,0,0\right )  ds \notag\\&= \frac 1 2 \frac 1 {x^3}+ \frac 1 {x^3} P\left(2+\epsilon  - \frac 1 x \right)- \frac 1 x \int_{0} ^{\frac 1 x} s P(2+\epsilon - s) ds. \notag
\end{align} Thus, we have this lemma.\qed

\vskip10pt

Lemma \label{Na} describes the graph of $P(x)$ as an application of Lemma \ref{Ga}. Thus, $P (x)$ and $\p_y R_{B_1} \left(x -1 -\frac \epsilon 2 ,0,0\right)$ is positive and decreasing for $x > 1$.

\begin{lem}\label{Na} The function $ P(x)$ satisfies $$ P (x) > 0, ~P'(x)<0, ~P''(x)>0, ~P'''(x)<0, ~P''''(x)>0$$ for $x>1$ and  $$\lim_{x \rightarrow \infty} P(x) = 0.$$
\end{lem}

\pf  To prove the decay of $P$ at infinity, we use the limits of  left- and right-hand sides of the equality \eqref{B} as $x$ approaches $\infty$. The interval $\left[1+\frac 1 2 ,2+\frac 1 2\right]$ contains a neighborhood of $2+\epsilon$, and the continuous $P(x)$ is bounded on the compact set $\left[1+\frac 1 2 ,2+\frac 1 2\right]$. Indeed \eqref{abcdefghi} means $\norm {P}_{L^{\infty}} \leq C \frac 1 {\epsilon}$ for some $C>0$.  Then, the equality \eqref{B} implies   $$\lim_{x \rightarrow \infty} P(x) = 0.$$

\par Next, we consider the positivity of $(-1)^n P^{(n)} (x)$ for $n=0,1,2,3,4$. Similarly with the precious lemma, we first study the properties of the reflection only with respect to a single inclusion. Second, we apply such properties to the case of two neighboring inclusions.

\par According to plan, we  consider the properties \eqref {lem2_2:A}, \eqref {lem2_2:B}, \eqref {lem2_2:C}, \eqref {lem2_2:D} and \eqref {lem2_2:E} of the reflection with respect to a single inclusion $B_0$ that denotes the unit sphere $B_1 (0,0,0)$. As defined before, $R_0 (h)$ denotes the reflection of a given harmonic function $h$ with respect to $B_0$. Suppose that  for $n=0,1,2,3,4,$ $$ \p_x^{n}\p_y h (x,0,0) \geq 0 \mbox{  on } (0,1].$$ We shall show that for $n=0,1,2,3,4,$ $$(-1)^n \p_x ^n \p_y R_0 ( h) (x,0,0)\geq 0 \mbox{  on } (0,1].$$ This means \eqref {lem2_2:A}, \eqref {lem2_2:B}, \eqref {lem2_2:C}, \eqref {lem2_2:D} and \eqref {lem2_2:E}. 

\par To do so, we use the equality \eqref{R_0_deriv} in the previous lemma or \eqref{R_0_deriv_1}. First, the positivity of $  \p_{y} R_0( h) (x, 0, 0)$ results from the equality immediately as  \begin{align}\p_{y} R_0( h ) (x, 0, 0) &= \frac 1 {x^3} \p_{y} h \left(\frac 1 x,0,0\right) - \frac 1 x \int_{0} ^{\frac 1 x  }  s \p_{y} h \left(s,0,0\right) ds \notag \\&= \frac 1 x \int_{0} ^{\frac 1 x  }  s  \left( 2 \p_{y} h \left(\frac 1 x,0,0\right) - \p_{y} h \left(s,0,0\right)\right) ds  \geq 0 \label{lem2_2:A} \end{align} due to the increasing property of $\p_{y} h \left(x,0,0\right) \geq 0 $. Second, dealing with  the decreasing property of $\p_{y} R_0 (h)(x, 0, 0)$, we take a derivative of  \eqref{R_0_deriv} and then, the increasing assumption of  $\p_{y} h(x, 0, 0)$  yields that \begin{align}
\p_{x}\p_{y} R_0( h) (x, 0, 0) &=   - \frac 2 {x^4}   \p_{y} h \left(\frac 1 x,0,0\right) + \frac 1 {x^2} \int_{0} ^{\frac 1 x  }  s \p_{y} h \left(s,0,0\right) ds - \frac 1 {x^5}  \p_x  \p_{y} h \left(\frac 1 x,0,0\right)  \notag \\&\leq - \frac 1 {x^2} \int_{0} ^{\frac 1 x  }  s  \left( 4 \p_{y} h \left(\frac 1 x,0,0\right) - \p_{y} h \left(s,0,0\right)\right) ds  \leq 0. \label{lem2_2:B}
\end{align} Thus, $\p_{y} R_0( h ) (x, 0, 0)$ is decreasing. Third, the concavity result can be also obtained in the same way. Thus,  \begin{align}
&\p_{x} ^2 \p_{y} R_0( h) (x, 0, 0) \notag \\ &=    \frac 7 {x^5}   \p_{y} h \left(\frac 1 x,0,0\right) - \frac 2 {x^3} \int_{0} ^{\frac 1 x  }  s \p_{y} h \left(s,0,0\right) ds + \frac 7 {x^6}  \p_x  \p_{y} h \left(\frac 1 x,0,0\right)    + \frac 1 {x^7}  \p_x ^2 \p_{y} h \left(\frac 1 x,0,0\right) \notag \\ &\geq \frac 1 {x^3} \int_{0} ^{\frac 1 x  }  s  \left( 14 \p_{y} h \left(\frac 1 x,0,0\right) - 2\p_{y} h \left(s,0,0\right)\right) ds  \geq 0.\label{lem2_2:C}
\end{align} Fourth, we have similarly \beq \p_{x} ^3 \p_{y} R_0( h) (x, 0, 0) \leq    -\frac {33} {x^6}   \p_{y} h \left(\frac 1 x,0,0\right) + \frac 6 {x^4} \int_{0} ^{\frac 1 x  }  s \p_{y} h \left(s,0,0\right) ds  \leq 0.\label{lem2_2:D}\eeq and \beq \p_{x} ^4 \p_{y} R_0( h) (x, 0, 0) \geq    \frac {192} {x^7}   \p_{y} h \left(\frac 1 x,0,0\right) - \frac {24} {x^5} \int_{0} ^{\frac 1 x  }  s \p_{y} h \left(s,0,0\right) ds  \geq 0.\label{lem2_2:E}\eeq

\par At last, we are ready to prove this lemma. By Lemma \ref{another_repre_sol}, 
  $$u(\Bx)= H (\Bx)+ R_1(H) (\Bx) + R_2(H) (\Bx) + R_2(R_1(H)) (\Bx) + R_1(R_2(H) )(\Bx) + \cdots , $$ and \begin{align*}\p_{y} u (x,0,0)=& \p_{y} H(x,0,0)+ \p_{y} R_1(H) (x,0,0) + \p_{y} R_2(H) (x,0,0) \notag\\&+ \p_{y} R_2(R_1(H)) (x,0,0)+ \p_{y} R_1(R_2(H) )(x,0,0) + \cdots,  \end{align*}
 where $R_1$ and $R_2$ are the reflections with respect to the insulated inclusions $B_1$ and $B_2$ as defined in \eqref{3-5_2015}. We apply  \eqref{lem2_2:A}, \eqref{lem2_2:B}, \eqref{lem2_2:C}, \eqref{lem2_2:D}, \eqref{lem2_2:E}.  Let $n=0,1,2,3,4$. Since $\p_y H(x,0,0) =1$, $$ (-1)^n  \p_{x} ^n \p_{y} R_1(H) (x,0,0) \geq 0 ~\mbox{ for } x > -\frac {\epsilon} 2 $$ and  $$ \p_{x} ^n \p_{y} R_2(H) (x,0,0) \geq 0 ~\mbox{ for } x < \frac {\epsilon} 2 .$$ In the same way, one can show by the mathematical induction that $$ (-1)^n  \p_{x} ^n \p_{y}( ( R_1R_2)^m (H)) (x,0,0) \geq 0, $$ $$ (-1)^n  \p_{x} ^n \p_{y}( R_1( R_2R_1)^m (H)) (x,0,0) \geq 0  $$ for $ x > -\frac {\epsilon} 2,$ and  $$ \p_{x} ^n \p_{y} ( (R_2 R_1)^m (H))) (x,0,0)\geq 0, $$ $$ \p_{x} ^n \p_{y} ( R_2 (R_1 R_2)^m (H))) (x,0,0) \geq 0  $$ for $ x < \frac {\epsilon} 2$, when $m\in \mathbb{N}$. Two representations \eqref {sol_decompos_2015} and \eqref {this_lemma_series} of the solution $u$ yield \begin{align*} \p_{y} R_{B_1} (x,0,0) = &\p_{y} R_1(H)(x,0,0) \\& + \p_{y} R_1(R_2(H)) (x,0,0) + \p_{y} R_1(R_2(R_1 (H))) (x,0,0)+\cdots,\end{align*} 
 and   by \eqref {prop32-proof-2015}, $$(-1)^n \p_x ^n \p_{y} R_{1} (H) \left(x ,0,0\right) >0.$$ We thus have $$(-1)^n \p_x ^n \p_{y} R_{B_1} (x,0,0) > 0 $$ for $ x > -\frac {\epsilon} 2$ that implies this lemma, since $P(x+ 1 +\frac {\epsilon} 2) = \p_{y} R_{B_1} \left(x ,0,0\right) $.
\qed

\vskip 10pt
\par  Another property of $P$ is provided by the following lemma based on the previous lemma. This property is used to reduce \eqref {B} into an ordinary differential equation in Lemma \ref{lem27}.
\begin{lem}\label{Da} For $n =1,2,3,4$,
$$\left| (x-1)^n P^{(n-1)} (x)\right| \lesssim (x-1) P\left(\frac {x-1}  2 +1 \right) $$ for any $x > 1$.
\end{lem}

\pf  For $n=1,2,3,4$, the decreasing property of  $|P^{(n-1)} |$ in $(1, \infty)$ is provided in the previous lemma.   In the case of $n=1$, it yields $$\left| (x-1) P (x)\right| \lesssim (x-1) P\left(\frac {x-1}  2 +1 \right) $$ for any $x > 1$.

Let $n$ be one of $2,3,4$.  By the mean value theorem, for any $x > 1$,  there exists $x_0 \in \left(x - \frac {x-1} {2^n} , x \right)$ such that 
$$ \left| P ^{(n-1)}(x_0)\right| = 2^{n-1}\left| \frac {P^{(n-2)} \left(x - \frac {x-1} {2^{n-1}} \right) - P^{(n-2)}  (x)} {x-1}\right| \lesssim \left| \frac {P^{(n-2)}  \left(x - \frac {x-1} {2^{n-1}}  \right)} {x-1}\right|,$$ since the value of $(-1)^{(n-2)} P^{(n-2)}$ is always  positive. It follows from the decreasing property of $|P^{(n-1)} |$ that
$$ \left| P ^{(n-1)}(x)\right| \leq  \left| P ^{(n-1)} (x_0)\right|  \lesssim \left| \frac {P^{(n-2)}\left(x - \frac {x-1} {2^{n-1}}\right)} {x-1}\right|.$$ When $n = 3$ or $4$, we continue this process so that
\begin{align*}\left| P ^{(n-1)}(x)\right| &\lesssim \left| \frac {P^{(n-2)}\left(x - \frac {x-1} {2^{n-1}}\right)} {x-1}\right|  \\& \lesssim \frac 1 {|x-1|} \left| \frac {P^{(n-3)}\left(x - \frac {x-1} {2^{n-2}}\right)- P^{(n-3)}\left(x - \frac {x-1} {2^{n-1}}\right) } {x-1}\right|\\& \lesssim \left| \frac {P^{(n-3)}\left(x - \frac {x-1} {2^{n-2}}\right)} {(x-1)^2}\right|     \\& \lesssim \cdots \lesssim \left| \frac {P\left(x - \frac {x-1} {2}\right) - P\left(x - \frac {x-1} {2^2}\right)} {(x-1)^{n-1}}\right|   \\&\lesssim \left| \frac {P\left(\frac {x-1}  2 +1 \right)} {(x-1)^{n-1}}\right|.\end{align*} Thus, we have this lemma. \qed

\vskip 10pt
The fundamental equation \eqref{B} can be rewritten as 
$$  \frac 1 2 \frac 1 {x^2}   -\int_{0} ^{\frac 1 x} s P(2+\epsilon - s) ds = xP(x) -  \frac 1 {x^2} P\left(2+\epsilon  - \frac 1 x \right)  .$$ The left-hand side is positive by the following lemma. The value of the left-hand side is very important, since the blow-up rate of $P$ is proportional to $$\frac 1 {\sqrt \epsilon}\left(\frac 1 2 \frac 1 {x^2}   -\int_{0} ^{\frac 1 x} s P(2+\epsilon - s) ds\right) \mbox{ at } x=1+ 2 \sqrt \epsilon.$$ Refer to Lemma \ref {337} for the details. 
\begin{lem} \label{integral}
\beq  \int_{0} ^{\frac 1 x} s P(2+\epsilon - s) ds <\frac 1 2 \frac 1 {x^2} \label{integral_inequal} \eeq
for any $x \in [1,2+\epsilon]$.
\end{lem}
\pf First, we show that the inequality \eqref {integral_inequal} is valid for $x$ on the restricted interval  $[1,p_2]$. Here, $p_2$ is the fixed point given in \eqref{p_1p_2}. Second, the inequality on $[1,2+\epsilon]$ is proved by contradiction.

\par According to plan, we prove that \beq  \int_{0} ^{\frac 1 x} s P(2+\epsilon - s) ds - \frac 1 2 \frac 1 {x^2} < 0 ~\mbox{ for any }
x \in [1,p_2].\label{integral_result1}\eeq It is easy to show that $x < 2+\epsilon  - \frac 1 x $ and $x\geq1$ for any $x \in [1,p_2)$. The decreasing property of $P$ in Lemma \ref{Na} yields $0> P\left(2+\epsilon  - \frac 1 x \right) -  P(x)$ and $P\left(2+\epsilon  - \frac 1 x \right) > 0 $ for any $x \in [1,p_2)$. By Lemma \ref {Ga}, \begin{align*}
 0>& x \left(\frac 1 {x^3} P\left(2+\epsilon  - \frac 1 x \right) -  P(x)\right)\\ =&\int_{0} ^{\frac 1 x} s P(2+\epsilon - s) ds - \frac 1 2 \frac 1 {x^2}
\end{align*} for any $x \in [1,p_2]$,  since $x\geq1$ and $2+\epsilon - \frac 1 {p_2} = p_2 >1$. Thus, we got the  result \eqref{integral_result1} restricted on $ [1,p_2]$.

\par Suppose that  $$\int_{0} ^{\frac 1 {x_0}} s P(2+\epsilon - s) ds - \frac 1 2 \frac 1 {x_0 ^2} \geq 0 $$ for some $x_0 \in (1,2+\epsilon]$.   By the mean value theorem, there exists a point $s_0 \in (0, \frac 1 {x_0})$ such that $P(2+\epsilon - s_0) \geq 1$, since $\int_{0} ^{\frac 1 {x_0}} s ds = \frac 1 2 \frac 1 {x_0 ^2}  $. The decreasing property $P$ yields
$$ P(2+\epsilon - s) \geq 1 $$  for any $s >  s_0$. For any $x \in [1,x_0]$, $\frac 1 x  > \frac 1 {x_0} > s_0$ so that  \begin{align*} & \int_{0} ^{\frac 1 x} s P(2+\epsilon - s) ds - \frac 1 2 \frac 1 {x^2} \\ &= \int_{\frac 1 {x_0}} ^{\frac 1 x} s P(2+\epsilon - s) ds - \frac 1 2 \left(\frac 1 {x^2}  - \frac 1 {x_0^2}  \right) + \int_{0} ^{\frac 1 {x_0}} s P(2+\epsilon - s) ds    - \frac 1 2  \frac 1 {x_0^2}   \\ &\geq \int_{\frac 1 {x_0}} ^{\frac 1 x} s ds - \frac 1 2 \left(\frac 1 {x^2}  - \frac 1 {x_0^2}  \right)  = 0. \end{align*}  This leads to a contradiction for the first result \eqref{integral_result1} in this proof. Thus, we have $$\int_{0} ^{\frac 1 {x}} s P(2+\epsilon - s) ds - \frac 1 2 \frac 1 {x ^2} < 0 $$ for any $x \in [1,2+\epsilon]$.
\qed
\vskip 15pt

\begin{rem} \label{integral_rem} It follows from Lemma \ref{integral} that 
$$  0 <\frac 1 2 - \int_{0} ^{\frac 1 x} s P(2+\epsilon - s) ds  $$
for any $x \in [1,2+\epsilon]$. Thus, we have   $$\int_{\epsilon} ^{\frac 1 2} P\left(1+s\right) ds \leq 1,$$ and the decreasing property of $P$ yields  $$\int_{\epsilon} ^{1} P \left(1+ s\right) ds \leq 3 $$ due to the smallness of  $\epsilon$.
\end{rem}

\par The function $P$ is a solution to a second-order ordinary differential equation as will be seen in Lemma \ref{lem27}. The solution $P$ is decomposed into a particular solution and a linear combination of two homogenous solutions. The following lemma is used to estimate the coefficients in the linear combination that are essential to find the blow-up rate of $P$. Refer to Proposition \ref{upper_bound_f} for the details.

\begin{lem} \label{prop_P}  
$$ 3 \gamma \sqrt {\epsilon} P\left(1+ \gamma \sqrt {\epsilon} - \gamma^2  {\epsilon} \right) + \left((\gamma^2 -1 )\epsilon - \gamma^3 \epsilon^{\frac 3 2}\right) P'( 1+ \gamma \sqrt {\epsilon}) > 0 $$ in $2 < \gamma < \frac 1 {10\sqrt {\epsilon}}$.
\end{lem}

\pf Applying Lemma \ref{integral} to \eqref{B}, we have 
\begin{align}
0&\leq P(x)- \frac 1 {x^3} P\left(2+\epsilon  - \frac 1 x \right) \notag\\& = \left( 1- \frac 1 {x^3}\right) P\left(2+\epsilon  - \frac 1 x \right) + \left( P(x)-  P\left(2+\epsilon  - \frac 1 x \right)\right). \label{prop_P_1}
\end{align} 

\par We estimate the ingredients in \eqref{prop_P_1}. Let $x= 1+ \gamma \sqrt \epsilon$, while $2 < \gamma  < \frac 1 {10\sqrt {\epsilon}}$. Then, $$1- \frac 1 {x^3} \leq 3 \gamma \sqrt {\epsilon},$$ since $ \gamma < \frac 1 {10\sqrt {\epsilon}}$.   Note that $  x > 2+\epsilon  - \frac 1 x$, since $x > p_2 $ due to the conditions $ \gamma > 2$ and \eqref{p_1p_2}.  The mean value theorem provides the existence of $x_0 \in \left( 2+\epsilon  - \frac 1 x, x \right) $ such that \begin{align*} P(x)-  P\left(2+\epsilon  - \frac 1 x \right) &= \left(   x + \frac 1 x - 2 - \epsilon\right)P'(x_0)\\&\leq \left(  x + \frac 1 x -2 - \epsilon \right)P'(x)\leq 0 ,\end{align*} since $P'(x_0) \leq P'(x) \leq 0 $ by the monotonic property of $P'$ in Lemma \ref{Na}.  Since $2< \gamma < \frac 1 {10\sqrt {\epsilon}} $, we also have $$ - \left(  x + \frac 1 x -2 - \epsilon \right) \leq  - \left((\gamma^2 -1 )\epsilon - \gamma^3 \epsilon^{\frac 3 2}\right) \leq 0 $$ and 
 $2+\epsilon  - \frac 1 x \geq  1+ \gamma \sqrt {\epsilon} - \gamma^2  {\epsilon}$ that implies $$0<  P\left(2+\epsilon  - \frac 1 x \right) \leq P \left(1+ \gamma \sqrt {\epsilon} - \gamma^2  {\epsilon} \right)$$ due to the decreasing property of $P$.
 Applying these bounds above to \eqref{prop_P_1}, we have this lemma. \qed

\vskip 15pt
 We consider the property of $ 2+ \epsilon - \frac 1 x$ in the equation \eqref{B}, since the equation is the key ingredient in the proof of the first main result.
\begin{lem} \label{prop_x_n}  Suppose that the sequence $\{x_n\}$ satisfies 
$$ x_1= 1$$ and $$ x_{n+1} = 2+ \epsilon - \frac 1 {x_n}~\mbox{for}~n\in \mathbb{N}.$$
Then, for any $ ~n\leq \frac 1 {2\sqrt {\epsilon}} $, $$  x_n = 1 + (n-1) \epsilon + o_n $$
and $$ |  o_n |\leq 10 n^2 \epsilon \sqrt {\epsilon}.$$
\end{lem}

\pf One can show that \beq   x_n =  {p_2}  +\frac {2+ \epsilon - 2 p_2 }  { c_0 d^{n-1}+1 } ,\label {lem:2_7_x_n}\eeq
where $$p_2= 1+ \sqrt {\epsilon +  \left(\frac {\epsilon}2\right)^2 } + \frac {\epsilon}2$$ as given in \eqref{p_1p_2},  $$ c_0 = \frac {1-p_2 +\epsilon } {1-p_2}  \mbox{ and } d= \frac {p_2} {2+\epsilon -{p_2}} .$$

\par We estimate $p_2$, $c_0$ and $d$  in \eqref{lem:2_7_x_n}.  Thus, $$p_2 = 1 + \sqrt {\epsilon} + \frac1 2 {\epsilon} +  \frac 1 8 \epsilon \sqrt {\epsilon}  + O(\epsilon^{2+\frac 1 2}),$$ $$c_0 = 1 - \sqrt {\epsilon} + \frac1 2 {\epsilon} + O(\epsilon \sqrt {\epsilon}) \geq 1 - \sqrt {\epsilon} + \frac1 2 {\epsilon} $$ and $$ d= 1+2 \sqrt {\epsilon} + 2 {\epsilon} + O(\epsilon \sqrt {\epsilon}) > 1+2 \sqrt {\epsilon} + 2 {\epsilon} $$   for small $\epsilon>0$.  It has been proved in \cite{LY} that 
$$  1+ (n-1)x   \leq (1+x)^{n-1} \leq 1 + (n-1) x + (n-1)^2 x^2 ,$$
supposed that $x\in (0,2)$ and $(1+x)^{n-1} \leq 2$. Since $n\leq \frac 1 {2\sqrt {\epsilon}} $, $c_0 d^{n-1}$ can be estimated as 
$$ 1 + (2n-3) \sqrt {\epsilon}  \leq  c_0 d^{n-1} \leq 1 + (2n-3) \sqrt {\epsilon} + 10 n^2 \epsilon .$$
Applying these bounds above to  \eqref{lem:2_7_x_n}, we have 
\begin{align*} x_n &= 1+\sqrt \epsilon + \frac 1 2 \epsilon - \sqrt{\epsilon} \left(1 + \left( \frac 3 2 - n \right)\sqrt {\epsilon}\right) +  o_n \\ &=1 + (n-1) \epsilon + o_n \end{align*} where $$ |  o_n |\leq 10 n^2 \epsilon \sqrt {\epsilon} . $$
\qed

\subsection{Proof of Propostion \ref{prop_2_2}  } \label{subsection3-2}
The decreasing property of  $\p_{y} R_{B_1} \left(x ,0,0\right)$ was presented in Lemma \ref{Na}. We first show the existence of a constant  $r_0 $ regardless of $\epsilon$ such that   $$P (1 + r_0 \sqrt {\epsilon}) \lesssim  \frac 1 { \epsilon ^{\frac {2 - \sqrt 2 } {2} }}.$$ Second, a relation between the values $P(1)$ and $P (1 + r_0 \sqrt {\epsilon})$ is established in \eqref{last_step}.  Then, we can prove that  $$P (1 ) \lesssim  \frac 1 { \epsilon ^{\frac {2 - \sqrt 2 } {2} }}.$$ This implies the first bound \eqref {prop41thefirst2015} in Proposition \ref{prop_2_2}.  The second bound  \eqref {prop41thefirst2015-2} is also presented by virtue of  the  positivity and decreasing property of  $\p_{y} R_{B_1} \left(x ,0,0\right)$ in Lemma \ref{Na}.

\subsubsection {Estimate for $P (1 + r_0 \sqrt {\epsilon})$ for a large $r_0 >0$}
For the sake of convenience, we begin by defining $t$ and $f$ as
\begin{equation*}
\quad \left\{
\begin{array}{ll}
\ds t= x-1, \\
\ds f(t) = P(1+t) = P(x)
\end{array}
\right. 
\end{equation*} for $t \geq 0 $. 

\par Speaking of the scheme, the function $f$, defined in $[0,\infty)$, is a solution to the ordinary differential equation in Lemma \ref{lem27}. The function $f$ can be decomposed into three functions  in \eqref{thedecomp_f} as follows:
$$  f = f_p + C_{\alpha} f_\alpha + C_{\beta} f_{\beta},$$
where $f_p$ is a particular solution, and $f_\alpha$ and $f_\beta$ are two homogeneous solutions satisfying \beq f_\alpha (t) \simeq \frac 1 {t^{2-\sqrt 2}}  \mbox { and } f_\beta (t) \simeq \frac 1 {t^{2+\sqrt 2}} \label{2015_may_14}\eeq for $t \geq 10 \sqrt {\epsilon}$. The boundedness of $f_p$ is provided in Lemma \ref{particular}, and the boundedness of  $C_\alpha$ and the smallness of $ C_{\beta}$ are also derived by Lemma \ref{prop_P2}. Hence, we can estimate $$ P (1 + r_0 \sqrt {\epsilon})= f(r_0 \sqrt {\epsilon})  \lesssim f_{\alpha} (r_0 \sqrt {\epsilon}) \simeq \frac 1 {{\epsilon}^{\frac {2-\sqrt 2}2}}  $$ in Proposition \ref{upper_bound_f}  and Remark \ref{rem_3_5}. This is the scheme to estimate for $P (1 + r_0 \sqrt {\epsilon})$ for a large $r_0 >0$.

\par In the following lemma, we establish the ordinary differential equation which $f$ satisfies.
\begin{lem}\label{lem27}
$$(t^2 - \epsilon) f'' (t) + 5 t f'(t) + 2 f(t) = - \frac 1 {(1+t)^3} + g (t) $$
and $$|g (t)| \lesssim \left|t f\left(\frac t 2 \right)\right|$$
for any $ t > 10 \sqrt \epsilon$.
\end{lem}

\pf By Lemma \ref{Ga}, 
 $$ \frac 1 2 \frac 1 {x^2} =- \frac 1 {x^2} P\left(2+\epsilon  - \frac 1 x \right)+  \int_{0} ^{\frac 1 x} s P(2+\epsilon - s) ds + x P(x)$$ for  any $x \geq 1$.
Taking derivative, we have 
\begin{align}  - \frac 1 {x^3}= &- \frac 1 {x^4} P' \left(2+\epsilon - \frac 1 x\right)  + x P'(x) + P(x)+\frac 1{x^3} P\left(2+\epsilon -\frac 1 x \right)  \label{deri_P} \\ =&- \frac 1 {x^4}\left( P' \left(2+\epsilon - \frac 1 x\right) -P' (x) \right) + \left(- \frac 1 {x^4} + x \right)P' \left( x\right)  \notag \\&+ \left(1+ \frac 1{x^3}  \right)P(x) + \frac 1{x^3}  \left(P\left(2+\epsilon -\frac 1 x \right)  -   P(x)\right)  .   \notag\end{align} Since $ 2+\epsilon - \frac 1 x =x + \left(2+\epsilon - \frac 1 {1+t} -(1+t)\right)  = x  +( \epsilon - t^2) + O (t^3)$ in $0<t <1$,  the mean value theorem yields \begin{align}
&(t^2 - \epsilon) f'' (t) + 5 t f'(t) + 2 f(t) \notag \\& = - \frac 1 {(1+t)^3} + (3t + O(t^2) )f(t) +(-1+3t + O(t^2) ) (\epsilon - t^2 + O(t^3)) f'(t_1)  \notag\\& ~~~+ O(t^2) f'(t_2)  + 4t (\epsilon +O( t^2)) f''(t)   + O \left((\epsilon + t^2 )^2\right) f'''(t_3) \notag \\& = - \frac 1 {(1+t)^3} +  O(t) f(t) + O(t^2) f'(t_1)  + O(t^2) f'(t_2) + O(t^3) f''(t) + O(t^4) f'''(t_3) \notag\\&   = - \frac 1 {(1+t)^3} + g(t) \label{details_of_g}
\end{align} for $t >\sqrt \epsilon$, where the points $t_1$, $t_2$, $t_3$  are located between  $t + \epsilon - t^2$ and $ t $, and they are depending on $t$. Lemma \ref{Da} means that $$\left| t^n f^{(n-1)} (t)\right| \lesssim t f \left(\frac t 2 \right)$$for $t >0$. Thus, $$|g (t)| \lesssim \left|t f\left(\frac t 2 \right)\right|$$
for any $ t \in \left(10 \sqrt \epsilon, \frac 1 2 \right)$.
\qed

\vskip 20pt

\par Now, we consider the solution to 
\begin{equation} (t^2 - \epsilon) f'' (t) + 5 t f'(t) + 2 f(t) = - \frac 1 {(1+t)^3} + g (t) ~\mbox{ for  }t \geq 10 \sqrt {\epsilon} . \notag\end{equation} 
We shall find three proper functions $f_p$, $f_\alpha$ and $f_\beta$ that satisfy 
\beq (t^2 - \epsilon) f_p'' (t) + 5 t f_p'(t) + 2 f_p(t) = - \frac 1 {(1+t)^3} + g (t) \label{eqn_particular} \eeq
and 
$$ (t^2 - \epsilon) f_i '' (t) + 5 t f_i '(t) + 2 f_i(t) = 0$$ for $i=\alpha,~\beta$.
The general solution is decomposed into the three functions as follows:
\beq f = f_p + C_{\alpha} f_\alpha + C_{\beta} f_{\beta}, \label{thedecomp_f}\eeq
where $f_{\alpha}$ and $f_{\beta}$ is homogeneous solutions defined as
\begin{align} f_{\alpha} = t^{-2+ \sqrt 2} +  t^{-2+ \sqrt 2}  \sum_{n=1} ^{\infty}  \left( \frac {\epsilon}{t^2}\right)^n  \prod_{k=1}^{n}   \frac {(2k-\sqrt 2 ) (2k+1 - \sqrt 2 ) } {2k (2k-2\sqrt 2)},\label{d_f_alpha}\\ f_{\beta} = t^{-2- \sqrt 2} +  t^{-2- \sqrt 2}  \sum_{n=1} ^{\infty}  \left( \frac {\epsilon}{t^2}\right)^n  \prod_{k=1}^{n}   \frac {(2k+\sqrt 2 ) (2k+1 + \sqrt 2 ) } {2k (2k+2\sqrt 2)}\label{d_f_beta} \end{align} for $t \geq 10 \sqrt {\epsilon}$.

\par The functions $f_{\alpha}$ and $f_{\beta}$ can be established by induction. To do so, we regard $f_{\alpha}$ and $f_{\beta}$ as the sums $\sum_{n=0} ^{\infty} f_{\alpha n}$ and $\sum_{n=0} ^{\infty} f _{\beta n}$ , where $ f_{\alpha 0} = t^{-2+ \sqrt 2}  $ and $ f_{\beta 0}=  t^{-2- \sqrt 2}$ are the solutions to $ t^2 f_{i0}'' (t) + 5 t f_{i0}'(t) + 2 f_{i0}(t) = 0 ,$ and  $f_{i n}$ is the solution to $ t^2 f_{in}'' (t) + 5 t f_{in}'(t) + 2 f_{in}(t) = \epsilon f''_{i (n-1)} $ for $i=\alpha,~\beta$, and $n=1,2,3,\cdots$.  The functions $f_{\alpha}$ and $f_{\beta}$ are defined well on $[10\sqrt {\epsilon}, \infty)$, because $$\left|\frac {f_{\alpha n }(t) } {f_{\alpha (n-1)} (t) } \right|= \left(\frac {\epsilon}{t^2}   \right) \left|\frac {(2n-\sqrt 2 ) (2n+1 - \sqrt 2 ) } {2n (2n-2\sqrt 2)}\right|\leq 4 \left(\frac {\epsilon}{t^2}   \right)   , $$ $$\left|\frac {f_{\beta n }(t) } {f_{\beta (n-1)} (t) } \right|= \left(\frac {\epsilon}{t^2}   \right) \frac {(2n+\sqrt 2 ) (2n+1 + \sqrt 2 ) } {2n (2n+2\sqrt 2)} \leq 4 \left(\frac {\epsilon}{t^2}   \right)  $$ and the variable $t \geq 10\sqrt {\epsilon}$. Moreover, we have 
$$f_{\alpha} \simeq   t^{-2+ \sqrt 2}~\mbox{ and }~f_{\beta} \simeq   t^{-2- \sqrt 2}. $$

Dealing with \eqref{thedecomp_f}, we consider the contribution of $f_p$ to $f$. The boundedness of $f_p$ is derived in the following lemma. 
\begin{lem} \label{particular} There are a particular solution $f_p$ to \eqref{eqn_particular} and a constant $C_0$ such that 
$$|f_p (t)| \lesssim 1  $$ and $$|f_p '(t)| \lesssim \frac 1 t  $$ for any $t > C_0 \sqrt \epsilon$.
\end{lem}
\pf 
We shall find the sequence of functions $\{f_{pn}\}$ satisfying 
\begin{align}
& t^2  f_{p0}'' (t) + 5 tf_{p 0}'(t) + 2 f_{p0}(t) = G(t), \label{relat_f_{p0}}\\  & t^2  f_{p n}'' (t) + 5 tf_{p n}'(t) + 2 f_{pn}(t) =  \epsilon f_{p(n-1)}'' (t)\label{relat_f_{pn}}
\end{align} for $n=1,2,3,\cdots$, where $$G(t) =  - \frac 1 {(1+t)^3} + g (t).$$  The sum $\sum_{n=0} ^{\infty} f_{pn}$ is the desirable function $f_p$.    Without any loss of generality, we  assume in this proof that \beq |g (t)| \leq  t f\left(\frac t 2 \right)  \label{sim_lem27}\eeq for any $ t > 10 \sqrt \epsilon$.  This is a simplication of the inequality $|g (t)| \lesssim  t f\left(\frac t 2 \right)$ in Lemma \ref{lem27} for convenience.

 Define  $f_{p0}(t)$ as
\beq f_{p0}(t) =\frac {1}{t^{2 +\sqrt 2 }} \int_{10 \sqrt \epsilon} ^t \frac {1}{w^{1 -2\sqrt 2 }} \int_{10 \sqrt \epsilon} ^w G(s) s^{1 -\sqrt 2 } ds dw .\label{def_f_p0}\eeq Then, $f_{p0}(t)$ is a solution to \eqref {relat_f_{p0}}. Moreover, we estimate $\left| f_{p0}^{(n)}(t) \right| $ for $n=0,1,2$. By Remark \ref{integral_rem}, $\int_{2 \epsilon} ^{1}  f\left(\frac s 2\right) ds \leq 2$. Lemma \ref{lem27} and \eqref{sim_lem27} yield
\begin{align*}| f_{p0}(t) |  & \leq \frac {1}{t^{2 +\sqrt 2 }} \int_{10 \sqrt \epsilon} ^t \frac {1}{w^{1 -2\sqrt 2 }} \int_{10 \sqrt \epsilon} ^w  s^{1 -\sqrt 2 }  + f\left(\frac s 2\right) s^{2 -\sqrt 2 } ds dw \\&  \leq \frac 1 2+ \frac {1}{t^{2 +\sqrt 2 }} \int_{10 \sqrt \epsilon} ^t    w^{1 +\sqrt 2 } \int_{10 \sqrt \epsilon} ^w  f\left(\frac s 2\right) ds dw \leq 2 \end{align*} in  $10 \sqrt \epsilon <t<1$.
Taking the derivative of \eqref{def_f_p0},  we can get similarly   $$ | f_{p0}'(t) | \leq \frac {11} t $$ and by \eqref{relat_f_{p0}} $$  | f_{p0}''(t) | \leq \frac 1 {t^2}(60 + |g(t)|) $$ in  $10 \sqrt \epsilon <t<1$.

\par We also define $f_{p n }(t)$ as  \beq f_{p n }(t) =\frac {1}{t^{2 +\sqrt 2 }} \int_{\sqrt {120\epsilon} } ^t \frac {1}{w^{1 -2\sqrt 2 }} \int_{\sqrt {120\epsilon} } ^w \epsilon f_{p(n-1)}'' s^{1 -\sqrt 2 } ds dw. \label {def_f_pn}\eeq Then, $f_{p n }(t)$ is the solution to \eqref{relat_f_{pn}}. In the same way, we can prove by mathematical induction and \eqref{sim_lem27} that  for any $n=1,2,\cdots$, \beq | f_{pn}(t) | \leq  \frac {2}{2^{n}}, ~| f_{pn}'(t) | \leq  \frac {11}{2^n } \frac 1 t \label{f_pn_induction1}\eeq  and \beq | \epsilon f_{pn}''(t) | \leq  \frac {1}{2^n}  \frac {\epsilon} {t^2}  (60+ |g(t)|) \leq  \frac {1}{2^{n+1}}   (1+ |g(t)|) \label{f_pn_induction2}\eeq that is the right-hand side of \eqref{relat_f_{pn}}, while $ \sqrt {120 \epsilon} \leq t < \frac 1 2$. Hence, the sum $\sum_{n=0} ^{\infty} f_{p0}$ is well defined  and is the desirable function $f_p$. \qed

\vskip 15pt

 \par We shall consider the contributions of $f_{\alpha}$ and $f_{\beta}$ to $f$ in Proposition \ref{upper_bound_f}, since the boundedness of $f_p$ was derived in the previous lemma. To do so, we need Lemmas \ref{prop_P2} and  and \ref{minor_one}.
\begin{lem} \label{prop_P2}  
$$ 3 \gamma \sqrt {\epsilon}~ f\left( \gamma \sqrt {\epsilon} - 2\gamma^2  {\epsilon} \right) + \left((\gamma^2 -1 )\epsilon - \gamma^3 \epsilon^{\frac 3 2}\right) f'( \gamma \sqrt {\epsilon}) > 0 $$
in  $2 < \gamma < \frac 1 {10\sqrt {\epsilon}}$.
\end{lem} The lemma above is a rewritten version of Lemma \ref{prop_P}, since $f(t) = P(1+t)$.
\vskip 10pt

\par We also need the following lemma to prove Proposition \ref{upper_bound_f}.
\begin{lem} \label{minor_one} Suppose that the constants $ \widetilde M>0$, $ \widetilde C_{\alpha}$, $ \widetilde C_{\beta}>0$ and $ \widetilde C_{0}>0$ satisfy  $$\widetilde M +  \widetilde C_{\alpha} \frac 1 {t^{2- \sqrt 2}} \geq \widetilde C_{\beta}  \frac 1 {t^{2+ \sqrt 2}} ~\mbox{ for any}~t \geq \widetilde C_{0} \sqrt \epsilon.$$ Then, $$\frac 1 2 \left(\widetilde M +  \widetilde C_{\alpha} \frac 1 {t^{2- \sqrt 2}}\right) \geq \widetilde C_{\beta} \frac 1 {t^{2+ \sqrt 2}} ~\mbox{ for any}~t \geq 2 \widetilde C_{0} \sqrt \epsilon.$$
\end{lem} \pf  For any $t \geq \widetilde C_{0} \sqrt \epsilon$, $$\widetilde M {t^{2- \sqrt 2}}  +  \widetilde C_{\alpha}\geq \widetilde C_{\beta} \frac 1 {t^{2 \sqrt 2}} ~.$$ 
Let $t = s ( \widetilde C_{0} \sqrt \epsilon)$. For any $s \geq 2^{\frac 1 {2\sqrt 2}}$, Then,$$\widetilde M {t^{2- \sqrt 2}}  +  \widetilde C_{\alpha} \geq \widetilde M {\left( \widetilde C_{0} \sqrt \epsilon\right)^{2- \sqrt 2}}  +  \widetilde C_{\alpha}\geq   \widetilde C_{\beta} \frac 1 {{\left( \widetilde C_{0} \sqrt \epsilon\right)}^{2 \sqrt 2}} \geq  2 \widetilde C_{\beta} \frac 1 {t^{2 \sqrt 2}}.$$\qed

Using Lemmas \ref{prop_P2} and \ref{minor_one}, we estimate $f(t)$ in $\left[ \gamma_0 \sqrt \epsilon , \frac 1 {10}  \right)$ in the following proposition, supposed that $\gamma_0$ is sufficiently large regardless of $\epsilon$.  We shall prove in Subsection \ref{688subsubsection}  that $P\left(1+\gamma_0 \sqrt \epsilon \right) = f\left(\gamma_0 \sqrt \epsilon \right)$ has the same blow-up rate as $P(1)$. In this respect, the estimate for $f\left(\gamma_0 \sqrt \epsilon \right)$ is meaningful. 

\begin{prop} \label{upper_bound_f} There is a constant $\gamma_0 $ independent of $\epsilon$ such that 
$$ f (t)  \lesssim f_{\alpha} (t) \simeq \frac 1 {t^{2-\sqrt 2}}$$  for $t \in \left[ \gamma_0 \sqrt \epsilon , \frac 1 {10}  \right) $.
\end{prop}
\pf  We consider the constants $C_{\alpha}$ and $C_{\beta}$ in  the decomposition of $f$ as \beq f = f_p + C_{\alpha} f_{\alpha} + C_{\beta} f_{\beta}\label{decamp}\eeq as given in \eqref{thedecomp_f}. We recall the behavior of $f_p$, $f_{\alpha}$ and $f_{\beta}$. The boundedness of  $\norm {f_p (t) }_{L^{\infty} \left([C_0 \sqrt \epsilon, \frac 1 2 ]\right)}$ and $\norm { t f_p'(t) }_{L^{\infty} \left([C_0 \sqrt \epsilon, \frac 1 2] \right)}$ was presented in Lemma \ref {particular}. From the definitions \eqref{d_f_alpha} and \eqref{d_f_beta}, 
\beq f_{\alpha} (t) = \frac 1 {t^{2-\sqrt 2}} \left (1+ O \left( \frac {\epsilon} {t^2}\right)\right)\mbox{ and } f_{\beta} (t) = \frac 1 {t^{2+\sqrt 2}} \left (1+ O \left( \frac {\epsilon} {t^2}\right)\right) \label{falphafbetaestimate433}\eeq for $t > C_0 \sqrt {\epsilon}$, supposed that $C_0$ is sufficiently large regardless of $\epsilon$. In this proof, we consider two cases when $C_{\beta} \leq 0$ and $C_{\beta} \geq 0$, respectively.

\par In the first case when $C_{\beta} \leq 0$, the decomposition \eqref{decamp} and the positivity of $f$ in Lemma \ref{Na} yield that $  f_p + C_{\alpha} f_{\alpha} \geq - C_{\beta} f_{\beta}  >0 $. By the boundedness of $f_p$, there exists $M_1 >0$ such that \beq \frac 1 2 M_1 \geq \norm {f_p (t) }_{L^{\infty} \left([C_0 \sqrt \epsilon, \frac 1 2 ]\right)} . \label{Cbetaleq0}\eeq Then, $$  M_1 + C_{\alpha} f_{\alpha}  (t)\geq -C_{\beta} f_{\beta} (t)$$ for  $t > C_0 \sqrt {\epsilon}$.  Note that 
$ f_{\alpha} \simeq \frac 1 {t^{2-\sqrt 2}} $ and  $f_{\beta} \simeq \frac 1 {t^{2+\sqrt 2}} $. By Lemma \ref{minor_one}, there is $r_0 >0 $ such that  \beq  \frac 1 2 \left(M_1 + C_{\alpha} f_{\alpha}  (t) \right) \geq -C_{\beta} f_{\beta} (t) \label{eq:2_9}\eeq for $ t \geq r_0 \sqrt \epsilon > C_0 \sqrt \epsilon$. Thus, \begin{align*} 1 &\gtrsim  \int_{r_0 \sqrt \epsilon} ^{\frac 1 2 } f + M_1  dt \\& \geq \int_{r_0 \sqrt \epsilon} ^{\frac 1 2 }   \frac 1 2  C_{\alpha}  f_{\alpha} + \left(\frac 1 2 M_1 + f_p\right) + \left(\frac 1 2 \left(M_1 + C_{\alpha}f_{\alpha}\right)+ C_{\beta} f_{\beta}\right)   dt    \\& \geq   C_{\alpha} \int_{r_0 \sqrt \epsilon} ^{\frac 1 2 }  \frac 1 2   f_{\alpha} dt \gtrsim C_{\alpha}, \end{align*} since $ 1 \gtrsim  \int_{r_0 \sqrt \epsilon} ^{\frac 1 2 } f  dt$ by Remark \ref{integral_rem}. Hence, we use \eqref{eq:2_9}  to get
$$ f (t)\lesssim 1 + C_{\alpha} f_{\alpha}(t)  \lesssim  f_{\alpha} (t) $$ for $ t \geq r_0 \sqrt \epsilon$.

\par In the second case when $C_{\beta} \geq 0$,  we have $ f \geq  f_p + C_{\alpha} f_{\alpha}.$ Then, it follows from Remark \ref{integral_rem} and Lemma \ref {particular} that   
$$C_{\alpha}  \int _{C_0 \sqrt \epsilon}^{1} f_{\alpha} dt \leq \int _{C_0 \sqrt \epsilon}^{1} f - f_p dt \lesssim 1.$$ Thus, since $\int _{C_0 \sqrt \epsilon}^{1} f_{\alpha} dt \simeq 1$, the constant $$ C_{\alpha}  \lesssim 1.$$ On the other hand, we estimate $C_{\beta}$. Note that $ 3- (2-\sqrt 2) >0  \mbox{ and }  3- (2+\sqrt 2) <0 .$  Applying  \eqref{decamp} and  \eqref{falphafbetaestimate433} to Lemma \ref{prop_P2}, we can find a positive constant $M_2$ regardless of $\epsilon$ so that
\beq M_2 + C_{\alpha} \frac 1 {t^{2-\sqrt 2}} \gtrsim  C_{\beta}\frac 1 {t^{2+\sqrt 2}} \label{m_c_alpha_1}\eeq
and 
\beq M_2 + C_{\alpha}f_{\alpha} (t) \gtrsim  C_{\beta}f_{\beta}(t)\label{m_c_alpha_2}\eeq in $C_0 \sqrt \epsilon < t \leq \frac 1 {10}$, where the boundedness of $f_p$ is also used to get $M_2$. Thus, we have 
$$ 0 \leq f (t) \lesssim M_2 + C_{\alpha} f_{\alpha} (t)  \lesssim  f_{\alpha} (t) $$ for $t \in \left[ C_0 \sqrt \epsilon , \frac 1 {10}  \right] $.

\par Hence, regardless of whether  $C_{\beta} \geq 0 $ or  $C_{\beta} \leq 0 $, $$ f (t)\lesssim   f_{\alpha} (t) $$ in $r_0 \sqrt \epsilon \leq t \leq \frac 1 {10}$.

\qed

\vskip15pt

\begin{rem}\label{rem_3_5} It follows from Proposition \ref{upper_bound_f} that 
$$P(1+t) \lesssim \frac 1 {t^{2-\sqrt 2}} $$  in $r_0 \sqrt \epsilon \leq t \leq \frac 1 {10}$. Hence, we have 
$$P(1+r_0 \sqrt \epsilon ) \lesssim \frac 1 {\epsilon^{\frac {2-\sqrt 2}2}}.$$
\end{rem}

\subsubsection {Estimate for $P (1)$}  \label{688subsubsection}
The estimate for $P (1)$ results from $ P(1+r_0 \sqrt \epsilon ) \lesssim \frac 1 {\epsilon^{\frac {2-\sqrt 2}2}}$ which was presented in Proposition \ref{upper_bound_f} and Remark \ref {rem_3_5}. We first prove Lemmas \ref {337} and \ref{subsubsection_221} to show $P \left(1+ \frac 1 {100} \sqrt \epsilon \right)\lesssim \frac 1 { \epsilon ^{\frac {2 - \sqrt 2 } {2} }} $ and second establish the chain \eqref{last_step} of correlations between $P(1)$ and $P \left(1+ \frac 1 {100} \sqrt \epsilon \right)$.  Hence, we can obtain the estimate for $P(1)$ as follows:
$$ P(1) \lesssim  \frac 1 {\epsilon^{\frac {2-\sqrt 2}2}}.$$

\begin{lem}\label{337}
\beq 0 < \frac 1 2 -  \int_{1+ 2\sqrt {\epsilon}} ^{2+\epsilon} (2+\epsilon - t) P(t) dt  \lesssim  \frac 1 { \epsilon ^{\frac {1 - \sqrt 2 } {2} }} . \label {P_p_2} \eeq
\end{lem}
\pf Let $x_* = 1+ \gamma_0 \sqrt \epsilon$ where $\gamma_0$ was given in Proposition \ref{upper_bound_f} and $r_0 > 3$. Since $x_* > 2+\epsilon  - \frac 1 {x_*} $, the decreasing property of $P$ in Lemma \ref{Na} implies $P(x_*) \leq P\left(2+\epsilon  - \frac 1 {x_*} \right) $. By Lemmas \ref{integral} and \ref{Ga},
$$0< \frac 1 2 \frac 1 {x_* ^2} - \int_{0} ^{\frac 1 {x_*} } s P(2+\epsilon - s) ds  = {x_*}  \left(P({x_*} ) - \frac 1 {x_* ^3} P\left(2+\epsilon  - \frac 1 {x_*}  \right) \right). $$ By Proposition \ref{upper_bound_f}  and Remark \ref{rem_3_5}, 
\begin{align*}0&< P(x_*) - \frac 1 {x_*^3} P\left(2+\epsilon  - \frac 1 {x_*} \right)  \\ &=  P(x_*) - P\left(2+\epsilon  - \frac 1 {x_*} \right)  + \left(1- \frac 1 {x_*^3} \right)P\left(2+\epsilon  - \frac 1 {x_*} \right)\\&\lesssim    \left(1- \frac 1 {x_*^3} \right)P\left(2+\epsilon  - \frac 1 {x_*} \right)\lesssim  \sqrt {\epsilon} \frac 1{\epsilon^{\frac {2-\sqrt 2} 2}}  = \frac 1{\epsilon^{\frac {1-\sqrt 2} 2}} .\end{align*} By Remark \ref{integral_rem} and the positivity of $P$  in Lemma \ref{Na}, 
\begin{align*}0&<  \frac 1 2 -  \int_{1+ 2\sqrt {\epsilon}} ^{2+\epsilon} (2+\epsilon - t) P(t) dt  \\&\leq  \frac 1 2 \frac 1 {x_* ^2} - \int_{2+\epsilon - \frac 1{x_*}} ^{2+\epsilon} (2+\epsilon - t) P(t) dt + r_0 \sqrt \epsilon \lesssim \frac 1{\epsilon^{\frac {1-\sqrt 2} 2}} , \end{align*} since $1+ 2\sqrt {\epsilon}   < 2+\epsilon - \frac 1{x_*}$ due to $r_0 > 3$. \qed

\begin{lem} \label{subsubsection_221}  
\beq  P \left(1+ \frac 1 {100} \sqrt \epsilon \right)\lesssim \frac 1 { \epsilon ^{\frac {2 - \sqrt 2 } {2} }}.\label{result_subsection_221}\eeq
\end{lem}
\pf  Applying $x_a= 1+ \frac 1 {100} \sqrt \epsilon $ to Lemma \ref{Ga}, we have 
$$    \frac 1 2  - {x_a}^2 \int_{0} ^{\frac 1 {x_a}} s P(2+\epsilon - s) ds =  {x_a}^3  P(x)  - P\left(2+\epsilon  - \frac 1 {x_a} \right), $$ and note that $$x_a < 2+\epsilon  - \frac 1 {x_a} < 1 + 2 \sqrt {\epsilon}.$$ It follows from \eqref{P_p_2}  that \begin{align}
 \frac 1 { \epsilon ^{\frac {1 - \sqrt 2 } {2} }}  &\gtrsim \frac 1 2 - {x_a}^2 \int_{0} ^{\frac 1 {x_a}} s P(2+\epsilon - s) ds \label{baram} \\
& = {x_a}^3  P(x_a) - P\left(2+\epsilon  - \frac 1 {x_a} \right)\notag \\
& \geq ( {x_a}^3 - 1 ) P(x_a). \notag
\end{align}  Then, $$  P \left(1+ \frac 1 {100} \sqrt \epsilon \right)\lesssim \frac 1 { \epsilon ^{\frac {2 - \sqrt 2 } {2} }}.$$
 \qed

 \vskip 10pt

\par To estimate $P \left(1 \right)$, we establish a correlation between $P(1)$ and  $P \left(1+ \frac 1 {100} \sqrt \epsilon \right)$. The correlation is  a long chain based on 
 $$ P\left(2+\epsilon  - \frac 1 x \right) =  x^3 P(x) + x^2  \int_{0} ^{\frac 1 x} s P(2+\epsilon - s) ds - \frac 1 2 $$ in Lemma \ref{Ga}.  In Lemma \ref {prop_x_n}, the sequence $\{x_n\}$ was defined as  
$$ x_1= 1~\mbox{and}~x_{n+1} = 2+ \epsilon - \frac 1 {x_n}~\mbox{for}~n\in \mathbb{N}.$$ Then,
 $$ P\left(x_{n+1} \right) =  x_n^3 P(x_n) + x_n ^2  \int_{0} ^{\frac 1 {x_n}} s P(2+\epsilon - s) ds - \frac 1 2 $$ for $n\in \mathbb{N}$.
 Let $n_0 = \left[ \frac 1 {20 \sqrt {\epsilon}}\right]$. Then,
  \begin{align} & P\left(x_{n_0 + 1 } \right) +\sum_{n=1} ^{n_0-1} x_{n+1}^3 \cdots \cdot  x_{n_0}^3  P\left(x_{n+1} \right) \notag\\& =  x_{n_0}^3 P(x_{n_0}) + x_{n_0 }^2  \int_{0} ^{\frac 1 {x_{n_0}}} s P(2+\epsilon - s) ds - \frac 1 2 \notag  \\&~~~+ \sum_{n=1} ^{n_0 -1 }   x_{n+1}^3 \cdots \cdot  x_{n_0}^3 \left( x_n^3 P(x_n) +   x_{n}^2  \int_{0} ^{\frac 1 {x_n}} s P(2+\epsilon - s) ds - \frac 1 2\right). \notag  \end{align} By cancellation, 
 \begin{align} 
 & x_{1}^3 \cdot x_{2}^3 \cdots \cdot  x_{n_0}^3  P(1) \notag   \\ &=  P\left(x_{n_0 + 1 } \right) + \left( \frac 1 2  -  x_{n_0 }^2  \int_{0} ^{\frac 1 {x_{n_0}}} s P(2+\epsilon - s) ds  \right) \notag \\&~~~+ \sum_{n=1} ^{n_0 -1 }   x_{n+1}^3 \cdots \cdot  x_{n_0}^3 \left( \frac 1 2 - x_{n}^2  \int_{0} ^{\frac 1 {x_n}} s P(2+\epsilon - s) ds \right). \label{last_step}   \end{align} It follows from  Lemma  \ref {prop_x_n}  that $ 1+ \frac 1 {100} \sqrt {\epsilon} \leq x_{n_0 + 1} $. Thus, by Lemma \ref{subsubsection_221} and the decreasing property of $P$ in Lemma \ref{Na},  $$ P\left(x_{n_0 + 1 } \right) \leq P \left(1+ \frac 1 {100} \sqrt {\epsilon} \right) \lesssim  \frac 1 { \epsilon ^{\frac {2 - \sqrt 2 } {2} }}.$$ We use  Lemma  \ref {prop_x_n} again so that $$ 1 \leq  x_{1}^3 \cdot x_{2}^3 \cdots \cdot  x_{n_0}^3  \lesssim 1,$$ and $ x_n \geq 1$ for $n\in \mathbb{N}$.  Note that $1+2\sqrt {\epsilon} > p_2 >x_{n}$ for $n=1,\cdots, n_0$, where $p_2$  is the fixed point in \eqref{p_1p_2}. By Lemma \ref{P_p_2}, 
$$\frac 1 { \epsilon ^{\frac {1 - \sqrt 2 } {2} }}  \geq \frac 1 2 - {x_n}^2 \int_{0} ^{\frac 1 {x_n}} s P(2+\epsilon - s) ds $$ for $n=1,\cdots, n_0$.  Therefore, the bound \eqref{last_step} can be reduced into the desirable bound \eqref {prop41thefirst2015} as follows:
$$ \p_y R_{B_1} \left(-\frac \epsilon 2 , 0, 0\right) = P(1) \lesssim  \frac 1 { \epsilon ^{\frac {2 - \sqrt 2 } {2} }}. $$

\par As mentioned earlier in Subsection \ref{subsection3-2}, the second bound  \eqref {prop41thefirst2015-2} is also obtained immediately by the  positivity and decreasing property of  $\p_{y} R_{B_1} \left(x ,0,0\right)$ in Lemma \ref{Na}. \qed

\section{ Proof of Theorem \ref{main_thm_lower}} \label{sec5}
We establish the lower bound in Theorem \ref{main_thm_lower} under the assumption that  $$H(x,y,z)= y \mbox{ in }\mathbb{R}^3.$$ The upper bound of $\nabla u$ and the properties presented in Section \ref{sec4} are used to derive the lower bound. The directional derivative  $\p_y u$ is decomposed as  $$ \p_y u (x,0,0)= 1 + \p_y R_{B_1} (x,0,0) + \p_y R_{B_2} (x,0,0)$$ for $|x|<\frac \epsilon 2$.  As defined in Subsection \ref{subsection3-2}, $$ f(t) = P(1+t)= \p_y R_{B_1} \left(t-\frac {\epsilon} 2,0,0\right)$$ for any $t >0$. 

\par In this proof, we shall prove the existence of $t_0 >0$ with \beq f(t_0) \gtrsim  \frac 1 {\epsilon^{\frac  {2-\sqrt 2 }2 }}. \label{sec5-the1steq}\eeq By Lemma \ref{Na}, $ \p_y R_{B_2} (x,0,0)= \p_y R_{B_1} (-x,0,0)>0$ for $x<\frac \epsilon 2$ and $\p_y R_{B_1} (x,0,0)$ is decreasing for $x>-\frac \epsilon 2$. Then, \eqref {sec5-the1steq} implies $$ \p_y u (x_0,0,0) \gtrsim  \frac 1 {\epsilon^{\frac  {2-\sqrt 2 }2 }}$$ for some $x_0 \in \left(-\frac \epsilon 2, \frac \epsilon 2\right)$. This is the desirable lower bound.

\par To show the existence of $t_0$, we need a negative particular solution $f_p$ to the equation \eqref{eqn_particular}.

\begin{lem}\label{lemma5-1-AAA} Let $g(t)$ and $G(t)$ be as given in \eqref{eqn_particular} and \eqref{relat_f_{p0}}. There exist a negative function $f_p$, three positive constants $M$, $s_0 > 100$ and $S_0< \frac 1 {10}$ such that  \beq  (t^2 - \epsilon) f_p'' (t) + 5 t f_p'(t) + 2 f_p(t) = - \frac 1 {(1+t)^3} + g (t) = G(t)  \label{lemma5-1-eqn111} \eeq
and $$ - M < f_p (t) <- \frac 1 {300}$$ for any $t \in (s_0 \sqrt \epsilon,  S_0)$.  Here, $M$, $s_0 $ and $S_0$ are the constants regardless of small $\epsilon >0$.
\end{lem} It is worth mentioning that the function $f_p$ in this lemma is slightly different from the one in Lemma \ref{particular}.

\vskip 10pt

\par \pf By Lemma \ref{lem27} and Proposition \ref{upper_bound_f}, there are constants $s_1$ and $S_1$ such that $$ - \frac 1 {(1+t)^3} + g (t)= G(t)  \leq - \frac 1 {4} $$ for any $t \in (s_1 \sqrt {\epsilon}, S_1) $ and $S_1 < \frac 1 {10}$.

\par As defined in Lemma \ref{particular}, we also consider $$ f_p (t)= \sum_{n=0} f_{pn} (t).$$ In the similar way to \eqref{def_f_p0} and \eqref{def_f_pn}, we define $f_{pn}(t)$ as
$$ f_{p0}(t) =\frac {1}{t^{2 +\sqrt 2 }} \int_{s_1 \sqrt \epsilon} ^t \frac {1}{w^{1 -2\sqrt 2 }} \int_{s_1 \sqrt \epsilon} ^w G(s) s^{1 -\sqrt 2 } ds dw  $$ and $$f_{p n }(t) =\frac {1}{t^{2 +\sqrt 2 }} \int_{s_1\sqrt {\epsilon} } ^t \frac {1}{w^{1 -2\sqrt 2 }} \int_{s_1 \sqrt {\epsilon} } ^w \epsilon f_{p(n-1)}'' s^{1 -\sqrt 2 } ds dw $$ for any $t \in (s_1 \sqrt {\epsilon}, S_1) $ and $n=1,2,3,\cdots$.  For any $t \in (2s_1 \sqrt {\epsilon}, S_1) $,$$ -M_1 < f_{p0}(t) \leq  - \frac 1 {128},$$ since $G(t) \leq -\frac 1 4$. Here, $M_1$ is the constant regardless of small $\epsilon$.

\par We can use the mathematical induction to get the analogues of  \eqref{f_pn_induction1} and \eqref{f_pn_induction2}. By virtue of $\frac {\epsilon} {t^2}$ in \eqref{f_pn_induction2}, there are constants $s_0 >100$ and $S_0<\frac 1{10}$ such that 
$$ |f_{pn} (t)|+ t |f_{pn}' (t)|+t^2 |f_{pn} ''(t)|  \leq \frac 1 {2^n }  \frac 1 {300 }   $$ for any $t \in (s_0 \sqrt \epsilon, S_0) \subset (s_1 \sqrt \epsilon, S_1)$ and $n=1,2,3,\cdots$. Hence,  $f_p$ is the solution to \eqref{lemma5-1-eqn111} and satisfies $$ -M < f_p (t) < -\frac 1 {300} $$ on the open interval $ (s_0 \sqrt \epsilon, S_0).$  Here, $M$ is the constant regardless of small $\epsilon$. 

\qed

\vskip 15pt

\begin{lem}\label{lemma5-2-AAA}  Let $s_0$ and $S_0$ be as given in the previous lemma. Then, $$\int_{s_0 \sqrt \epsilon} ^{S_0}  f(t) dt \gtrsim 1 .$$
\end{lem}

\pf Lemma \ref {337} yields
$$ \frac 1 4 \leq  \int_{1+ 2\sqrt {\epsilon}} ^{2+\epsilon} (2+\epsilon - t) P(t) dt$$
for small $\epsilon>0$, since $P(x) \lesssim \frac 1 {\epsilon^{\frac {2-\sqrt 2 } 2}}$. Thus, $$ \frac 1 {10} \leq  \int_{s_0 \sqrt \epsilon} ^{1+\epsilon} f (t) dt$$ for small $\epsilon>0$. The decreasing property of $f$ in Lemma \ref{Na} yields  $$ 1 \simeq \frac 1 {10} \left( \frac{S_0 - s_0 \sqrt \epsilon} {1+ \epsilon - s_0 \sqrt \epsilon}\right)  \leq  \int_{s_0 \sqrt \epsilon} ^{S_0} f (t) dt .$$ Thus, we have this lemma. \qed

\vskip 15pt
\par Now, we recall the decomposition \eqref{thedecomp_f} as 
\beq f = f_p + C_{\alpha} f_\alpha + C_{\beta} f_{\beta}\label{f_beta_M_bound_ABC}\eeq where $f_{\alpha}$ and $f_{\beta}$ are defined as \eqref{d_f_alpha} and \eqref {d_f_beta}.  We consider two cases when $C_{\beta} < 0$ and when $C_{\beta} > 0$, separately.

\par In the first case when $C_{\beta} < 0$, the negativity of $f_p$ in  Lemma  \ref{lemma5-1-AAA} and Lemma \ref{lemma5-2-AAA} yield $$ C_{\alpha} \int_{s_0 \sqrt \epsilon} ^{S_0} f_{\alpha} dt \geq \int_{s_0 \sqrt \epsilon} ^{S_0}  C_{\alpha} f_\alpha + \left(f_p  + C_{\beta} f_{\beta}\right) dt  =  \int_{s_0 \sqrt \epsilon} ^{S_0} f dt \gtrsim 1$$ due to the positivity of $f_\beta$ in \eqref {d_f_beta} and \eqref {falphafbetaestimate433}. Here, $s_0$ and $S_0$ are the constants in Lemma \ref{lemma5-1-AAA}.  From the definition \eqref {d_f_alpha}, $\int_{s_0 \sqrt \epsilon} ^{S_0} f_{\alpha} dt \lesssim 1$ so that \beq C_{\alpha} \gtrsim 1 . \label{section5eqn54asdf}\eeq  On the other hand,  the positivity of $f$ in  Lemma \ref{Na} and the negativity of  $f_p$ in Lemma \ref{lemma5-1-AAA} yield $C_{\alpha} f_{\alpha}(t) \geq - C_{\beta} f_{\beta}(t) - f_p (t) \geq - C_{\beta} f_{\beta} (t) >0  $  on the  interval $ (s_0 \sqrt \epsilon, S_0).$  Then, $$- \frac {C_{\alpha} f_{\alpha}(2s_0 \sqrt \epsilon)} {C_{\beta} f_{\beta} (2s_0 \sqrt \epsilon)}  \geq 1 ,$$ and considering the defintions \eqref{d_f_alpha} and \eqref{d_f_beta},  $$ - \frac {C_{\alpha} } {C_{\beta} } {(2s_0 \sqrt \epsilon)^{2\sqrt 2}}  > \frac 1 4, $$ since the constant $s_0 > 100$.
Thus, we use the definitions of $f_\alpha$ and $f_\beta$ again to get
\beq  - \frac {C_{\alpha} f_{\alpha}(t)} {C_{\beta} f_{\beta} (t)}  \geq  -\frac 1 4  \frac {C_{\alpha} } {C_{\beta} } {t^{2\sqrt 2}} \geq - \frac 1 4  \frac {C_{\alpha} } {C_{\beta} } {(2s_0 \sqrt \epsilon)^{2\sqrt 2}} \left(\frac t {2s_0 \sqrt \epsilon} \right) ^{2\sqrt 2} > 2 \label{please_final} \eeq  for $t \geq 8 s_0 \sqrt \epsilon$. Applying  \eqref{please_final} to \eqref{f_beta_M_bound_ABC}, the inequality \eqref{section5eqn54asdf} implies
$$ f(8s_0 \sqrt \epsilon) \geq \frac 1 2 C_{\alpha} f_{\alpha}(8s_0 \sqrt \epsilon) -M \gtrsim  \frac 1 { \epsilon^{\frac {2-\sqrt 2} 2}}, $$ when $C_{\beta} < 0$. Here, the constant $M$ is given in  Lemma  \ref{lemma5-1-AAA}.

\par In the second case when $C_\beta \geq 0$, we shall prove that \beq C_\alpha \gtrsim 1. \label{section5eqn54asdf-2}\eeq
 By \eqref{m_c_alpha_1} and \eqref{m_c_alpha_2} in the proof of Proposition \ref{upper_bound_f}, there are positive constants $C_A$, $s_2$ and $S_2$ such that  \beq C_A \left( M + C_{\alpha} \frac 1 {t^{2-\sqrt 2}} \right)\geq  C_{\beta}\frac 1 {t^{2+\sqrt 2}} \label{dealing_f_beta}\eeq and $$ C_A \left(  M + C_{\alpha}f_{\alpha} (t) \right) \geq  C_{\beta}f_{\beta}(t)$$ for any $ t \in (s_2 \sqrt \epsilon, S_2) \subset (s_0 \sqrt \epsilon  , S_0) $,  while the constant $M$ was given in Lemma \ref{lemma5-1-AAA}. To prove \eqref{section5eqn54asdf-2}, we subdivide this second case into two subcases, due to  \eqref{dealing_f_beta}.  First, if $$C_A M \geq \frac 1 2 C_\beta  \frac 1 {(2 s_2 \sqrt \epsilon)^{2+\sqrt 2}},$$ then, the definition of $f_\beta$  in \eqref {d_f_beta} yields
$$C_{\beta}\int_{2 s_{2} \sqrt \epsilon} ^{S_2} f_{\beta}  dt\lesssim \sqrt \epsilon.$$ In the same way as  Lemma \ref{lemma5-2-AAA}, the negativity of $f_p$ in Lemma \ref{lemma5-1-AAA}  and \eqref {f_beta_M_bound_ABC} yield $$1 \lesssim \int_{2 s_2 \sqrt \epsilon} ^{S_2} f dt  \lesssim C_{\alpha}\int_{2 s_2 \sqrt \epsilon} ^{S_2} f_{\alpha} dt + \sqrt \epsilon  .$$ Thus, since $ \int_{2 s_2 \sqrt \epsilon} ^{S_2} f_{\alpha} dt \lesssim 1 $, we have 
$$ C_{\alpha} \gtrsim 1.$$   Second, otherwise, if $$C_A M \leq \frac 1 2 C_\beta  \frac 1 {(2 s_2 \sqrt \epsilon)^{2+\sqrt 2}},$$ then \eqref{dealing_f_beta} implies $2C_A C_\alpha \frac 1 {(2 s_2 \sqrt \epsilon)^{2-\sqrt 2}}  \geq C_\beta \frac 1 {(2 s_2 \sqrt \epsilon)^{2+\sqrt 2}}.$ For any $t >  2 s_{2} \sqrt \epsilon$  $$C_\alpha \frac 1 {t^{2-\sqrt 2}}  \gtrsim C_\beta \frac 1 {t^{2+\sqrt 2}}.$$ Thus,  $$ C_\alpha \int_{2 s_{2} \sqrt \epsilon} ^{S_2} f_{\alpha}  dt   \gtrsim \int_{2 s_{2} \sqrt \epsilon} ^{S_2} C_{\alpha}f_{\alpha}+ C_{\beta}f_{\beta}  dt \geq \int_{2 s_2 \sqrt \epsilon} ^{S_2} f dt   \gtrsim 1 ,$$  since $f_p < -\frac 1 {300}$. This means $$C_{\alpha} \gtrsim 1.$$  Hence, it follows from  \eqref{f_beta_M_bound_ABC} that  
$$ f( 2 s_{2} \sqrt \epsilon) \geq C_{\alpha} f_{\alpha}( 2 s_{2} \sqrt \epsilon) - M \gtrsim  \frac 1 { \epsilon^{\frac {2-\sqrt 2} 2}}$$ when $C_{\beta} \geq 0$.

\par Therefore, as a result, there is a point $t_0$ with $$f(t_0) \gtrsim  \frac 1 {\epsilon^{\frac  {2-\sqrt 2 }2 }}$$ regardless of whether  $C_\beta \geq 0$ or $C_\beta < 0$. This is the desirable \eqref{sec5-the1steq}. \qed

\subsection*{Acknowledgement}
The author is deeply grateful to Professor H. Kang for his help and discussion, and also appreciates the regard of Professor M. Lim and her student S. Yu. After the completion of the first draft of this paper, Lim and Yu informed the author that they might be able to derive the same blow-up rate result. It is worthy mentioning that all results in this paper were obtained independently of their works.

\end{document}